\documentclass[12pt]{article}
\usepackage{tikz}
\usepackage{color}
\usepackage{amsmath}
\usepackage{amssymb}
\usepackage{enumitem}
\usepackage{graphicx}

\def\D{\mathrm{D}}
\def\d{\mathrm{d}}
\def\tw{\mathrm{tw}}
\def\red{\mathrm{red}}
\def\val{\mathrm{val}}
\def\log{\mathrm{log}}

\def\vir{\mathrm {vir}}

\def\virt{\mathrm{vir}}
\def\P{\mathsf{P}}
\def\PP{\mathbb{P}}
\def\DR{\mathsf{DR}}

\def\CP{{{\mathbb {CP}}}}

\def\ooC{{\mathcal{C}}}
\def\cL{{\mathcal{L}}}

\def\cO{\mathcal{O}}
\def\oM{\overline{\mathcal{M}}}
\def\cM{{\mathcal{M}}}
\def\C{{\mathcal{C}}}

\def\Z{\mathbb{Z}}

\def\ZZ{\mathcal{Z}}

\def\C{\mathbb{C}}
\def\Q{\mathbb{Q}}
\def\qed{{\hfill $\Diamond$}}
\def\cS{{\mathcal S}}
\def\cg{{\mathcal G}}

\def\Aut{{\rm Aut}}

\def\E{\mathrm{E}}
\def\n{\mathrm{n}}
\def\L{\mathrm{L}}
\def\V{\mathrm{V}}
\def\H{\mathrm{H}}
\def\g{\mathrm{g}}
\def\G{\mathsf{G}}
\def\rarr{\rightarrow}
\def\D{\mathsf{D}}

\def\nice{\displaystyle}
\def\tC{\widetilde{\mathsf{C}}}
\def\ccoarseC{{\mathsf{C}}}

\def\coarseL{\mathsf{L}}
\def\tL{\widetilde{{\mathcal L}}}

\def\ooMM{{\mathfrak{M}}}

\def\ooCC{{\mathfrak{C}}}

\def\ch{{\rm ch}}
\def\td{{\rm td}}

\usepackage{theorem}

\newtheorem{theorem}{Theorem}
\newtheorem{prop}[theorem]{Proposition}
\newtheorem{proposition}[theorem]{Proposition}
\newtheorem{corollary}[theorem]{Corollary}
\newtheorem{lemma}[theorem]{Lemma}

{\theorembodyfont{\rmfamily}
\newtheorem{definition}{Definition}

}

\title{Double ramification cycles with target varieties}
\author{F. Janda, R. Pandharipande, A. Pixton, D. Zvonkine}
\date{August 2020}

\begin{document}

\maketitle

\vspace{-20pt}

\begin{abstract} 
Let $X$ be a nonsingular projective algebraic
variety over $\mathbb{C}$, and let 
$\oM_{g,n,\beta}(X)$ be the moduli space of 
stable maps $$f:(C,x_1,\ldots,x_n) \to X$$ 
from genus~$g$, $n$-pointed curves $C$ to $X$ of degree $\beta$.
Let $S$ be a line bundle on $X$.
Let $A = (a_1, \dots, a_n)$ be a vector of integers
which satisfy 
$$\sum_{i=1}^n  a_i = \int_\beta c_1(S)\, .$$ 
Consider  the following  condition: 
 {\em the line bundle $f^*S$  has a meromorphic section with zeroes and poles 
exactly at the marked points $x_i$  with 
orders prescribed by the integers $a_i$}. In other words, we require  $f^*S \left( -\sum_{i=1}^{n} a_i x_i \right)$ to be the trivial line bundle on $C$.

A compactification of the space of maps based upon the above condition 
is given by the moduli space of stable maps  to {\em rubber} over $X$ 
and is denoted by $\oM^\sim_{g,A,\beta}(X,S)$. 
The moduli space carries a virtual fundamental class 
$$[\oM^\sim_{g,A,\beta}(X,S)]^\virt\in A_*\left(
\oM^\sim_{g,A,\beta}(X,S)\right)$$ in Gromov-Witten theory. The main result of the paper
is an explicit formula (in 
tautological classes) for the push-forward via the forgetful morphism of 
$[\oM^\sim_{g,A,\beta}(X,S)]^\virt$ to $\oM_{g,n,\beta}(X)$. 
In case $X$ is a point, the result here specializes to 
Pixton's formula for the double ramification cycle proven in \cite{JPPZ}.
Several applications of the new formula are given.

\end{abstract}

\parskip=5pt
\baselineskip=15pt
\pagebreak

\vspace{-20pt}

\setcounter{tocdepth}{1} 
\tableofcontents

%%%%%%%%%%%%%%%%%%%%%%%%%%%%%%%%%%%%%%%%%%%%%%%%%%%%%%%%%%%%%%%%%%%%%%%%%%%

%\baselineskip{12pt}

\setcounter{section}{-1}
\section{Introduction} \label{Sec:intro}

\subsection{Double ramification cycles}
\label{Ssec:DRclassic}

Let $A = (a_1, \dots, a_n)$ be a vector of $n$ integers satisfying  $$\sum_{i=1}^n a_i = 0\,. $$ In the moduli space $\cM_{g,n}$ of 
nonsingular curves of genus $g$ with $n$ marked points,
 consider the substack defined by the following classical condition: 
\begin{equation}\label{g445}
\left\{ (C, x_1, \dots, x_n) \subset \cM_{g,n} \; \rule[-0.8em]{0.05em}{2em} \; \cO_C\Big(\sum_{i=1}^n a_i x_i\Big) \simeq \cO_C \right\}\, . 
\end{equation}
%nd it is natural to ask what is the homology class of its closure or of its natural compactification in $\oM_{g,n}$.
From the point of view of relative Gromov-Witten theory, 
the most natural compactification of the substack \eqref{g445} 
is the space $\oM^{\sim}_{g,A}$ of stable maps to
{\em rubber} \cite{GrV,JLi}: stable maps to
 $\CP^1$ relative to $0$ and $\infty$ modulo the $\C^*$ action on $\CP^1$.

 The rubber moduli space carries a natural virtual fundamental class $\left[\oM^{\sim}_{g,A}\right]^\virt$ of dimension $2g-3+n$. The push-forward 
via the canonical morphism
$$ \epsilon:\oM^{\sim}_{g,A} \rightarrow  \oM_{g,n}$$
is the {\em double ramification cycle} 
\begin{equation}\label{relggww}
\epsilon_*\left[\oM^{\sim}_{g,A}\right]^\virt\, =\, \mathsf{DR}_{g,A}\ \in A_{2g-3+n}(\oM_{g,n})
\, .
\end{equation}
The double ramfication cycle
$\mathsf{DR}_{g,A}$
can also be defined via log stable maps (and was
motivated in part by 
Symplectic Field Theory \cite{EGH}). 
%All different approaches all involve virtual fundamental classes
%of map spaces, and all three definitions agree.

The classical approach to the locus \eqref{g445} is via Abel-Jacobi theory
for the universal curve. However, extending the Abel-Jacobi map over 
the boundary $\oM_{g,n} \setminus \cM_{g,n}$ of the moduli space of curves is
{\em not} straightforward.  Approaches by Marcus-Wise \cite{MW} and
Holmes \cite{Hol} (motivated by log geometry), nevertheless, provide a partial
resolution of the Abel-Jacobi which is sufficient to define a double
ramification cycle. The result also agrees with  definition \eqref{relggww}.

Eliashberg posed the question of computing $\mathsf{DR}_{g,n}$ in 2001. The hope for a possible formula
was strengthened 
in  \cite{FP} where the double ramification cycle
was proven to lie in the
tautogical ring of $\oM_{g,n}$.
Calculations on
the open set of compact type curves were given in \cite{GruZak,Hain}.
Integrals against the double ramification cycle{\footnote{Termed {\em rubber integrals} in the papers \cite{OP1,OP2,OP3}.}}
played a fundamental role in the
solution of the Gromov-Witten theory for target curves \cite{OP1,OP2,OP3}.

A complete formula for  $\mathsf{DR}_{g,A}$ in the
tautological ring of $\oM_{g,n}$
was conjectured by
Pixton in 2014 and proven in \cite{JPPZ} via 
Gromov-Witten theory. Pixton's formula expresses $\mathsf{DR}_{g,A}$
directly as a sum over stable graphs $\Gamma$  indexing the
boundary strata of $\oM_{g,n}$. The contribution of each stable
graph $\Gamma$ is the constant term  of a polynomial in $r$ naturally associated
to the combinatorics of $\Gamma$ and $A$. The proof of \cite{JPPZ}
was obtained by studying the  Gromov-Witten theory of
the target $\mathbb{P}^1$ with an orbifold $B\mathbb{Z}_r$-point at $0\in \mathbb{P}^1$ and
a relative point at $\infty\in \mathbb{P}^1$ in the  $r\rightarrow \infty$
limit.

Pixton's formula's opened new directions in the subject: new
formulas for Hodge classes \cite[Section 3]{JPPZ}, new relations  in the tautological ring
of $\oM_{g,n}$  \cite{cj}, new connections to the loci of meromorphic differentials
\cite[Appendix]{FarP},
and connections to new integrable hierarchies \cite{BRS}.
For a sampling of the subsequent 
study and applications, see \cite{BHPSS,B,BGR, BSSZ,
 Cav, CGJZ,FanWu,HKP,HPS,ObPix,Pix18,Sch}. We refer the
reader to \cite[Section 0]{JPPZ} and \cite[Section 5]{RPSLC} for more
leisurely introductions to the subject.

The double ramification cycle study above concerns 
the Gromov-Witten theory of a {\em point}.{\footnote{The moduli
space of stable maps to a point is $\oM_{g,n}$.}}
Our goal
here is to develop a full theory of double ramification
cycles for general nonsingular projective target varieties $X$.

Let 
$\oM_{g,n,\beta}(X)$ be the moduli space of 
stable maps $$f:(C,x_1,\ldots,x_n) \to X$$ 
from genus~$g$, $n$-pointed curves $C$ to $X$ of degree $\beta$.
Let $S$ be a line bundle on $X$.
Let $A = (a_1, \dots, a_n)$ be a vector of integers
which satisfy 
$$\sum_{i=1}^n  a_i = \int_\beta c_1(S)\, .$$ 
Consider  the following  condition analogous to \eqref{g445}: 
 {\em the line bundle $f^*S$  has a meromorphic section with zeroes and poles 
exactly at the marked points $x_i$  with 
orders prescribed by the integers $a_i$}. In other words, we require  $f^*S \left( -\sum_{i=1}^{n} a_i x_i \right)$ to be the trivial line bundle on $C$. 
Rubber maps with target $X$ provide a natural compactification of the
locus of solutions and define an $X$-valued double ramification
cycle in $A_*(\oM_{g,n,\beta}(X))$.
Our main result is a complete formula 
for the $X$-valued double ramification cycle  which generalizes
the structure of Pixton's formula.

The formal definition of the $X$-valued double ramification
cycle is given in Section \ref{Ssec:RubberToX}. After a discussion of
$X$-valued stable graphs and tautological classes in Sections
\ref{xvalsg} and \ref{Ssec:PsiXiEta}, our formula for the $X$-valued double
ramification cycles is presented in Section \ref{Ssec:MainFormula}.
In case $X$ is a
point, we recover our previous study \cite{JPPZ}.

The
double ramification cycle construction for target varieties
plays a crucial role in relative Gromov-Witten theory. 
Since the answer for the $X$-valued double ramification
cycles takes such a simple form, new directions are again opened
in the subject:
\begin{enumerate}
\item[$\bullet$] The study of the tautological ring of
  $\oM_{g,n,\beta}(X)$ is suggested by a theory of relations \cite{Bae}
  parallel to the case of point \cite{cj}.

\item[$\bullet$] For a pair $(X,D)$ where $D$ is nonsingular
  divisor, the relationship between the relative Gromov-Witten theory
  and the orbifold Gromov-Witten theory of the root stack is beautifully
  settled in \cite{Tseng1,Tseng2}.
  The study of the $D$-valued double ramification
  cycle here plays a crucial role. 
  
\item[$\bullet$]
The  $\mathbb{CP}^N$-valued $\mathsf{DR}$-cycle in the limit
$N\rightarrow \infty$, suitably interpreted, is 
a universal $\mathsf{DR}$-cycle on the moduli space of line
bundles on curves. The universal Abel-Jacobi theory on the Picard stack
\cite{BHPSS}
is both motivated by and dependent upon our calculation of $\mathbb{CP}^N$-valued $\mathsf{DR}$-cycles.
\end{enumerate}
Of course, the $X$-valued formula also leads immediately to
simple derivations of older results in Gromov-Witten theory.
Applications are discussed at the end of the paper in Section \ref{appp}.

%In a previous paper~\cite{JPPZ} we proved a conjecture by the second author providing an explicit formula for the DR-cycle in terms of tautological classes. Here we extend the notion of DR-cycles to stable maps with a smooth target manifold~$X$ and prove a similar formula for this cycle.

\subsection{Rubber maps with target~$X$}
\label{Ssec:RubberToX}
Let $X$ be a nonsingular projective variety over $\C$. 
Let 
$ S \rightarrow X $
be a line bundle, and let 
$$\PP(\cO_X \oplus S) \rightarrow X$$
be the canonically associated $\CP^1$-bundle over $X$.
Let 
$$ D_0\, ,\, D_\infty \, \subset \, \PP(\cO_X \oplus S)$$
be the divisors defined by the projectivizations of
the loci $\cO_X \oplus \{0\}$ and $\{0 \} \oplus S$
respectively. We will call $D_0$ the $0$-divisor and 
 $D_\infty$ the $\infty$-divisor.

Let $C$ be a nonsingular curve with $n$ marked points,  and let
$$f : C \to X$$ be an algebraic map of degree $\beta \in H_2(X,\Z)$. 
Furthermore, let $s$ be a nonzero meromorphic section of $f^*S$ over $C$, defined up to a multiplicative constant, with zeros and poles belonging
 to the set of marked points of $C$. We denote the orders of zeros and poles by $$A = (a_1, \dots, a_n)\, .$$ If the $i$th marking
 is neither a zero nor a pole, we set $a_i=0$. We have
 $$\sum_{i=1}^n a_i = \int_\beta c_1(S)\,. $$

The pair $(f,s)$ defines a map to  {\em  rubber with target $X$}. 
Let 
$$\oM^{\sim}_{g,A,\beta}(X,S)$$  be
the compact moduli space of stable maps to 
{rubber with target $X$}. 
A general stable map to rubber with target~$X$ 
is a map to a rubber chain of $\CP^1$-bundles $\PP(\C \oplus S)$ over~$X$ 
attached along their $0$- and $\infty$-divisors. 
The space of rubber maps with target~$X$ 
carries a perfect obstruction theory and a virtual fundamental class, see~\cite{JLi,Li2,MauPan} for a detailed discussion. 
Let  $$\epsilon: \oM^{\sim}_{g,A,\beta}(X,S) \to \oM_{g,n,\beta}(X)$$ 
be the morphism obtained by the projection of the rubber  map to $X$ (and the
contraction 
of the resulting unstable components).

%The {\em double ramification data} in the case of a target manifold~$X$ consists of a list of $n$ integers $A = (a_1, \dots, a_n)$ satisfying $\sum a_i  = c_1(S) \cdot \beta$ and, in addition, a degree $\beta \in H_2(X,\Z)$.

\begin{definition} \label{Def:main}
The {\em $X$-valued double ramification cycle} is 
the $\epsilon$ push-forward of the virtual fundamental class
of the moduli space of stable maps to rubber over $X$:
$$
\DR_{g,A,\beta}(X,S) \, =\, \epsilon_*
\left[ \oM^{\sim}_{g,A,\beta}(X,S) \right]^\vir\ 
 \in A_{\text{vdim}(g,n,\beta)-g}(\oM_{g,n,\beta}(X))\, .$$
\end{definition}
The virtual dimension
 $\text{vdim}(g,n,\beta)$  of
$\oM_{g,n,\beta}(X)$ is determined by
$$\text{vdim}(g,n,\beta) = \int_\beta c_1(X) + (g-1)(\text{dim}_\C(X) -3)
+n\, .$$

%\subsection{Virtual substacks}
%\label{Ssec:VirtualSubstacks}
%Our formula for the $X$-valued DR-cycle will be in terms of {\em $X$-valued stable graphs}, where each stable graph actually denotes a {\em virtual Cartier substack} (a generalization of Jun Li's C-divisors from~\cite{Li2}, Definition 3.1.) of $\oM_{g,n,\beta}(X)$. A {\em C-divisor} is simply a line bundle $E$ with a section; it is actually enough to define the line bundle in the neighborhood of the zero locus of the section, since it can always be extended by the trivial line bundle elsewhere. 

%If a stack $B$ is endowed with a perfect obstruction theory and a virtual class $[B]^\vir$, then one can easily define a perfect obstruction theory and a virtual class on the zero locus of the section of the C-divisor. We have then $[D]^\vir = [B]^\vir \cdot c_1(E)$. 

%Similarly, in higher codimension, given a vector bundle $E$ with a section, the zero locus $Z \subset B$ of the section is endowed with a natural virtual fundamental class $[B]^\vir \cdot c_{\rm top}(E)$ and called a {\em virtual substack}. The main advantage of this notion is that intersections and self-intersections of virtual substacks are naturally defined.

\subsection{$X$-valued stable graphs}\label{xvalsg}
We define the set $\G_{g,n,\beta}(X)$ of {\em $X$-valued stable graphs} as
follows. A graph $\Gamma \in \G_{g,n,\beta}(X)$ consists of the data
$$\Gamma=(\V\, ,\ \H\, ,\ \L\, , \ \mathrm{g}:\V \rarr \Z_{\geq 0}\, ,
\ v:\H\rarr \V\, , 
\ \iota : \H\rarr \H,\ \beta: \V \to H_2(X,\Z))$$
satisfying the properties:
\begin{enumerate}
\item[(i)] $\V$ is a vertex set with a genus function $\g:\V\to \Z_{\geq 0}$,
\item[(ii)] $\H$ is a half-edge set equipped with a 
vertex assignment $v:\H \to \V$ and an involution $\iota$,
\item[(iii)] $\E$, the edge set, is defined by the
2-cycles of $\iota$ in $\H$ (self-edges at vertices
are permitted),
\item[(iv)] $\L$, the set of legs, is defined by the fixed points of $\iota$ and is
placed in bijective correspondence with a set of $n$ markings,
\item[(v)] the pair $(\V,\E)$ defines a {\em connected} graph
satisfying the genus condition 
$$\nice \sum_{v \in \V} \g(v) + h^1(\Gamma) = g\, ,$$
\item[(vi)] for each vertex $v$, the stability condition holds: if $\beta(v)=0$,
then
$$2\g(v)-2+ \n(v) >0\, ,$$
where $\n(v)$ is the valence of $\Gamma$ at $v$ including 
both edges and legs,
\item[(vii)] the
degree condition holds:
$$\nice \sum_{v \in \V} \beta(v) = \beta\, .$$
\end{enumerate}
To emphasize $\Gamma$, the notation $\V(\Gamma)$, $\H(\Gamma)$, $\L(\Gamma)$, and $\E(\Gamma)$ will also be to used for the vertex, half-edges, legs, and
edges of $\Gamma$.

An automorphism of $\Gamma \in \G_{g,n,\beta}(X)$ 
consists of automorphisms
of the sets $\V$ and $\H$ which leave invariant the
structures $\L$, $\mathrm{g}$, $v$, $\iota$, and $\beta$.
Let $\Aut(\Gamma)$ denote the automorphism group of $\Gamma$.

An $X$-valued stable graph $\Gamma$ determines a moduli space $\oM_\Gamma$ of stable maps with the degenerations forced by the graph together with a 
canonical map, 
$$j_\Gamma : \oM_\Gamma \to \oM_{g,n,\beta}(X)\, .$$
%of  degree $|\Aut(\Gamma)|$. 
The moduli space $\oM_\Gamma$ is the substack of 
the product
$$
\oM_\Gamma\subset \prod_{v \in \V} \oM_{\g(v), \n(v), \beta(v)}(X)$$
cut out by the inverse images of the diagonal 
$\Delta \subset X \times X$
under the evaluations maps associated to the edges $e=(h,h')\in \E$,
$$
\prod_{v \in \V} \oM_{\g(v), \n(v), \beta(v)}(X)\ \stackrel{\text{ev}_e\ }{\longrightarrow}\
X \times X\, .$$
The moduli space $\oM_\Gamma$ carries a natural virtual fundamental class $\left[ \oM_\Gamma \right]^\virt$ defined by the refined intersection, 
\begin{equation}\label{pp334}
 [\oM_\Gamma]^\virt \, =\, 
\prod_{e\in \E} \text{ev}_e^{-1}(\Delta) \ \cap \ 
\prod_{v \in \V} \left[\oM_{\g(v), \n(v), \beta(v)}(X) \right]^\virt\, .
\end{equation}
%For every edge $e = (h', h'')$ of $\Gamma$, we multiply \eqref{pp334} by 
%the pull-back of the Poincar\'e dual cohomology class of the diagonal $\Delta \subset X \times X$ under the evaluation map $\ev_{h'} \times \ev_{h''}$. This product defines the class $\left[ \oM_\Gamma \right]^\virt$. 

%It is important to note that $j_* \oM_\Gamma$ is a virtual substack of $\oM_{g,n,\beta}(X)$. Indeed,  each edge $e$ of the graph defines a line bundle $N_e = T' \otimes T''$, where $T'$ and $T''$ are the tangent lines to the branches. The sum $E = \oplus_{e \in \E} N_e$ defines the deformation vector bundle. The class $j_* \left[ \oM_\Gamma \right]^\virt$ coincides with the virtual substack class. In particular,  intersections of substacks determined by prestable graphs are well-defined. 

\subsection{Tautological $\psi$, $\xi$, and $\eta$ classes} \label{Ssec:PsiXiEta}

The universal curve $$\pi:\ooC_{g,n,\beta}(X) \to \oM_{g,n,\beta}(X)$$ carries two natural line bundles: the relative dualizing sheaf $\omega_{\pi}$ 
and the pull-back $f^*S$  of the line bundle $S$ via the
universal map,
$$f: \ooC_{g,n,\beta}(X) \rightarrow X\, .$$
%For notationally simplicity, we will call the pull-back $S$ instead
%of $f^*(S)$. 

Let $s_i$ be the $i$th section of the universal curve, let 
$$D_i\subset \ooC_{g,n,\beta}(X)$$
be  the corresponding divisor, and let
$$\omega_\log = \omega_\pi\Big(\sum_{i=1}^n D_i\Big)$$ be the relative logarithmic line bundle 
with first Chern class $c_1 (\omega_\log)$. 
Let $$\xi = c_1(f^*S)$$ be the first Chern class of the pull-back of $S$.

\begin{definition} \label{Not:PsiXiEta}
The following classes in $\oM_{g,n,\beta}(X)$ are obtained from the
universal curve
$\ooC_{g,n,\beta}(X)$:
\begin{itemize}
    \item $\psi_i = c_1(s_i^* \omega_\pi)\, $,
    \item $\xi_i = c_1(s_i^* f^*S)\, $,
    \item $\eta_{a,b} = \pi_*\left(c_1(\omega_\log)^a \xi^b\right)\, $.
\end{itemize}
All can be viewed as either cohomology classes or as operational
Chow classes in $A^*(\oM_{g,n,\beta}(X))$. We will use the term
{\em Chow cohomology} for operational Chow.
\end{definition}

Since the class $\eta_{0,2} = \pi_*(\xi^2)$ will play a prominent role in our formulas for the $X$-valued $\mathsf{DR}$-cycle, we will use the
notational convention
$$\eta= \pi_*(\xi^2)\, .$$
The standard $\kappa$ classes are defined by the $\pi$ push-forwards
of powers of $c_1(\omega_\log)$, so we have
$$\eta_{a,0} = \kappa_{a-1}\, .$$

%DZ: The following paragraph is added:
Consider the moduli $\oM_\Gamma$ of stable maps described by a stable graph~$\Gamma$. Let $e$ be an edge of the graph composed of two half-edges $h_1$ and $h_2$. The space $\oM_\Gamma$ carries two natural cohomology classes $\psi_{h_1}$ and $\psi_{h_2}$ (associated to the two cotangent lines at
the node corresponding to the edge $e$) and a cohomology class $\xi_e = c_1(s_e^*S)$, where $s_e$ is the section of the universal curve determined by the same node.

\begin{definition}
A {\em decorated $X$-valued stable graph} $[\Gamma,\gamma]$ 
%of type $(g,n,\beta)$ 
is an 
$X$-valued stable graph $\Gamma \in \G_{g,n,\beta}(X)$ together with the following
decoration data $\gamma$:
\begin{itemize}
\item each leg $i\in \L$ is decorated with a monomial $\psi_i^a \xi_i^b$,
\item each half-edge $h\in \H\setminus \L$ is 
decorated
%{\footnote{The assignment of $\xi$ classes to half-edges in
%$\H\setminus \L$
%is certainly possible, but such classes are not needed in our formulas.}}
with a monomial $\psi^a_h$,
\item each edge $e\in \E$ is decorated with a monomial $\xi^a_e$,
\item each vertex in $\V$ is decorated with a monomial in the variables 
$\{\eta_{a,b}\}_{a+b \geq 2}$.
\end{itemize}
%We denote the set of decorated $X$-valued stable graphs by $\D_{g,n,\beta}(X)$.
\end{definition}

Let $\D\G_{g,n,\beta}(X)$ be the set of decorated $X$-valued stable graphs.
To each decorated graph $$[\Gamma,\gamma]\in \D\G_{g,n,\beta}(X)\,,$$ we assign the cycle class $j_{\Gamma*}[\gamma]$
obtained
via the push-forward via
$$j_\Gamma : \oM_\Gamma \rightarrow \oM_{g,n,\beta}(X)$$
of
 the action of the product of the
 $\psi$, $\xi$, and $\eta$ decorations
 on  $\left[ \oM_\Gamma \right]^\virt$,
\begin{equation}\label{dd66}
j_{\Gamma*}[\gamma] \, \stackrel{\text{def}}{=}\, j_{\Gamma*} \left( \gamma \cap \left[ \oM_\Gamma \right]^\virt\right) 
\, \in\,  A_*(\oM_{g,n,\beta}(X))\, .
\end{equation}
Our formula for the $X$-valued $\mathsf{DR}$-cycle 
is a sum of  cycle classes  \eqref{dd66}
assigned
to decorated $X$-valued stable graphs.

\subsection{Weightings mod $r$} \label{Ssec:weightings}
Following the notation of  Sections \ref{Ssec:RubberToX}-\ref{Ssec:PsiXiEta}, let 
$S \to X$ be a line bundle on a nonsingular projective variety
$X$. Fix the data 
$$g \geq 0\, , \ \ \beta \in H_2(X, \Z)\, ,\ \ A = (a_1, \dots, a_n)$$
subject to the condition 
$$\sum_{i=1}^n a_i =
\int_{\beta} c_1(S)\, .$$
Let $\Gamma \in \G_{g,n,\beta}(X)$ be an $X$-valued stable graph, and
let $r$ be a positive integer.

\begin{definition} A {\em weighting mod~$r$} of $\Gamma$ is a function on the set of half-edges,
$$ w:\H(\Gamma) \rightarrow \{0,1, \ldots, r-1\}\, ,$$
which satisfies the following three properties:

\begin{enumerate}
\item[(i)] $\forall i\in \L(\Gamma)$, corresponding to
 the marking $i\in \{1,\ldots, n\}$,
$$w(i)=a_i  \mod r\, ,$$
\item[(ii)] $\forall e \in \E(\Gamma)$, corresponding to two half-edges
$h,h' \in \H(\Gamma)$,
$$w(h)+w(h')=0 \mod r\, ,$$
\item[(iii)] $\forall v\in \V(\Gamma)$,
$$\sum_{v(h)= v} w(h)= \int_{\beta(v)} c_1(S) \mod r\, ,$$ 
where the sum is taken over {\em all} $\mathsf{n}(v)$ half-edges incident 
to $v$.\end{enumerate}
\end{definition}

%\begin{enumerate}
%\item[(i)] $\forall h_i\in \L(\Gamma)$ we have $w(h_i) = a_i \pmod r$,
%$$
%\begin{array}{cll}
%w(h_i)  &= a_i \pmod r &\text{ for } 1 %\leq i \leq n, \\
%w(h_i) & = 0  &\text{ for } n+1 \leq i %\leq n+m;
%\end{array}
%$$
%\item[(ii)] $\forall e \in \E(\Gamma)$, corresponding to two half-edges
%$h,h' \in \H(\Gamma)$,
%$$w(h)+w(h')=0 \pmod r\, ,$$
%\item[(iii)] $\forall v\in \V(\Gamma)$,
%$$\sum_{h \in \nu^{-1}(v)} w(h)= (\beta_v , \xi) \pmod r\, ,$$ 
%\end{enumerate}
%where the sum is taken over {\em all} $\mathsf{n}(v)$ half-edges incident to~$v$.
%\end{definition}

We denote by $\mathsf{W}_{\Gamma,r}$ the finite set of all possible weightings
mod $r$
of
$\Gamma$. The set $\mathsf{W}_{\Gamma,r}$ has cardinality $r^{h^1(\Gamma)}$.
We view $r$ as a {\em regularization parameter}.

\subsection{The double ramification formula} \label{Ssec:MainFormula}
We denote by
$\P_{g,A,\beta}^{d,r}(X,S)\in A_{\text{vdim}(g,n,\beta)-d}(\oM_{g,n,\beta}(X))$
the degree $d$ component of the tautological class 
\begin{multline*}
\hspace{-10pt}\sum_{
\substack{\Gamma\in \G_{g,n,\beta}(X) \\
w\in \mathsf{W}_{\Gamma,r}}
}
\frac{r^{-h^1(\Gamma)}}{|\Aut(\Gamma)| }
\;
j_{\Gamma*}\Bigg[
\prod_{i=1}^n \exp\left(\frac12 a_i^2 \psi_i + a_i \xi_i \right)
\prod_{v \in \V(\Gamma)} \exp\left(-\frac12 \eta(v) \right)
\\ \hspace{+10pt}
\prod_{e=(h,h')\in \E(\Gamma)}
\frac{1-\exp\left(-\frac{w(h)w(h')}2(\psi_h+\psi_{h'})\right)}{\psi_h + \psi_{h'}} \Bigg]\, .
\end{multline*} 
Inside the push-forward in the above formula, the first product 
$$\prod_{i=1}^n \exp\left(\frac{1}{2}a_i^2 \psi_{h_i}+a_i \xi_{h_i}\right)\, $$
is over $h\in \L(\Gamma)$
via the correspondence of legs and markings.
The class $\eta(v)$ is the $\eta_{0,2}$ class 
of Definition \ref{Not:PsiXiEta} associated to the vertex.
The third product is over all $e\in \E(\Gamma)$.
The factor 
$$\frac{1-\exp\left(-\frac{w(h)w(h')}{2}(\psi_h+\psi_{h'})\right)}{\psi_h + \psi_{h'}}$$
is well-defined since 
\begin{enumerate}
\item[$\bullet$]
the denominator formally divides
the numerator,
\item[$\bullet$] the factor is symmetric in $h$ and $h'$.
\end{enumerate}
No edge orientation is necessary.

The following fundamental polynomiality property of $\P_{g,A,\beta}^{d,r}(X,S)$ 
 is
parallel to Pixton's polynomiality in \cite[Appendix]{JPPZ} and
is a consequence of \cite[Proposition $3''$]{JPPZ}.

\begin{prop} For fixed $g$, $A$, $\beta$, and $d$, the \label{pply}
class
$$\P_{g,A,\beta}^{d,r}(X,S) \in A_*(\oM_{g,n,\beta}(X))$$
is polynomial in $r$ (for all sufficiently large $r$).
\end{prop}

We denote by $\P_{g,A,\beta}^d(X,S)$ the value at $r=0$ 
of the polynomial associated to $\P_{g,A,\beta}^{d,r}(X,S)$ by Proposition~\ref{pply}. In other words, $\P_{g,A,\beta}^d(X,S)$ is the {\em constant} term of the associated polynomial in $r$. 

The main result of the paper is a 
 formula for the $X$-valued double ramification cycle parallel{\footnote{Our
handling of the prefactor $2^{-g}$ in \cite[Theorem 1]{JPPZ} differs here.
The factors of $2$ are now placed in the definition of $\P_{g,A,\beta}^{d,r}$.}} to
Pixton's proposal in case $X$ is a point.

\begin{theorem} \label{Thm:main}
Let $X$ be a nonsingular projective variety with line
bundle $$S\rightarrow X\, .$$
For $g\geq 0$, double ramification data $A$, and $\beta \in H_2(X, \Z)$, 
we have
$$\DR_{g,A,\beta}(X,S) =  \P_{g,A,\beta}^g(X,S)\, \in \, A_{\text{\em vdim}(g,n,\beta)-g}(\oM_{g,n,\beta}(X))\, .$$ 
\end{theorem}

\subsection{Strategy of proof}

Consider the projective bundle $\PP(\cO_X \oplus S) \to X$. 
By applying Cadman's $r$th root construction \cite{Cadman} 
to the $0$-divisor $D_0$, we obtain a bundle 
\begin{equation}\label{kk99}
\PP(X,S)[r]\rightarrow X
\end{equation} 
where every fiber is a projective line with a single stacky point 
of stabilizer $\Z/r \Z$. 

Our proof of Theorem \ref{Thm:main} 
is obtained by studying the Gromov-Witten theory of the bundle \eqref{kk99}
relative to the $\infty$-divisor $D_\infty$. 
The virtual localization formula for the orbifold/relative
geometry $$\PP(X,S)[r]/D_\infty$$ yields 
relations for every $r$ which depend polynomially on $r$ (for all sufficiently
large $r$). 
After setting $r=0$, we obtain the equality of Theorem \ref{Thm:main}.
The argument has two main parts.

\vspace{5pt}
\noindent (i) Let $\oM^r_{g,A,\beta}(X,S)$ be the moduli space of stable maps
 $$f:C \to X$$ endowed with an $r${th} tensor root $L$ of the line bundle 
$f^*S (-\sum_{i=1}^n a_i x_i)$. Furthermore, let 
$$\pi : \ooC^r_{g,A,\beta}(S) \to \oM^r_{g,A,\beta}(X,S)$$ 
be the universal curve, and let $\cL$ be the universal $r$th root over the universal curve. A crucial step is to prove that
the push-forward of $$r\cdot c_g(-R\pi_*\cL)$$ to $\oM_{g,n,\beta}(X)$ 
is a polynomial in $r$ (for all sufficiently large $r$) and
that the polynomial has the {\em same} constant term as the 
polynomial $\P_{g,A,\beta}^{g,r}(X,S)$. 
Our formula for the $X$-valued $\mathsf{DR}$-cycle therefore
has a geometric interpretation in terms of the top Chern class $c_g(-R\pi_*\cL)$. 

Contrary to the case of ordinary $\mathsf{DR}$-cycles studied in \cite{JPPZ}, 
for the case of $X$-valued $\mathsf{DR}$-cycles, we cannot 
use Chiodo's formulas \cite{Chiodo2} to deduce polynomiality. Instead,
 we adapt Chiodo's computations to our geometric setting
 in Sections \ref{Ssec:GRR} and \ref{Ssec:GRRpoly}. 

\vspace{5pt}
\noindent (ii)
 We use the localization formula \cite{GrP} for the virtual
fundamental 
class of the moduli space of stable maps to the orbifold/relative geometry
$$\PP(X,S)[r]/D_\infty\, .$$ 
The positive (respectively, negative) coefficients $a_i$
 specify the monodromy conditions over the $0$-divisor (respectively, 
the tangency conditions along the $\infty$-divisor). 

The moduli space $\oM^r_{g,A,\beta}(X,S)$ appears in the
 localization formula. Indeed, the space of stable maps to the $\C^*$-invariant locus corresponding 
to the stacky $0$-divisor $D_0$ is precisely $\oM^r_{g,A,\beta}(X,S)$. 
The push-forward of the localization formula to $\oM_{g,n,\beta}(X)$ is a Laurent series in the equivariant parameter $t$ and in $r$. The coefficient of $t^{-1}r^0$  must vanish by geometric considerations. We prove that the relation obtained from
the coefficient of $t^{-1}r^0$  has {\em only two terms}: 
\begin{itemize}
\item The first is the
constant term in $r$ of the push-forward of $r\cdot c_g(-R\pi_*\cL)$ to $\oM_{g,n,\beta}(X)$. 
\item
The second term is the $X$-valued double ramification cycle 
$\DR_{g,A,\beta}(X,S)$ with a minus sign. 
\end{itemize}
The vanishing of
the sum of the two terms yields Theorem~\ref{Thm:main}. 

\subsection{Notation table}
To help the reader, we list here the symbols used for the various
spaces which arise in the paper:
\begin{itemize}
    \item $\PP(X,S)$ the $\CP^1$-bundle $\PP(\cO_X \oplus S)$ over~$X$,
    \item $\PP(X,S)[r]$ is the outcome of applying Cadman's $r$th root construction to the 0-divisor $D_0 \subset \PP(X,S)$,
    \item $\oM_{g,n,\beta}(X)$ is the space of stable maps $f:(C,x_1,\ldots,x_n) \to X$,
    \item $\oM^r_{g,A,\beta}(X,S)$ is the space of stable maps $f:(C,x_1,\ldots,x_n) \to X$ together with an $r$th root of $f^*S(-\sum_{i=1}^n a_ix_i)$,
     \item $\ooMM_{g,n}$ is the stack of prestable curves,
    \item $\ooMM^r_{g,n}$ is the stack of twisted prestable curves,
    \item $\ooMM_{g,n}^Z$ is the stack of prestable curves together with a degree~0 line bundle~$Z$,
    \item $\ooMM_{g,n}^{r,Z,{\text{triv}}}$ is the stack of prestable twisted curves together with a degree~0 line bundle~$Z$ where the stabilizer of every point of the twisted curve acts trivially in the fibers of the line bundle,

    \item $\ooMM_{g,n}^{r,L}$ is the stack of prestable twisted curves together with a degree~0 line bundle $L$ with no conditions on the stabilizers.
\end{itemize}

\subsection{Acknowledgments} 
We are grateful to Y.~Bae, A.~Buryak, R.~Cavalieri, A.~Chiodo, E.~Clader, 
H. Fan,
G.~Farkas, J.~Gu\'er\'e, A.~Kresch, G.~Oberdieck, D.~Oprea, J.~Schmitt, H.-H. Tseng,
and L. Wu for discussions about double ramification cycles and related topics. The {\em Workshop on Pixton's conjectures} at ETH Z\"urich  in October 2014 played in an important role in our collaboration. F.~J., A.~P., and D.~Z. have been frequent guests of the Forschungsinstitut f\"ur Mathematik (FIM) at ETHZ.

%F.~J. was partially supported by the Swiss National Science Foundation grant SNF-200021143274. 
R.~P. was partially supported by SNF-200020-162928,
SNF-200020-182181, ERC-2017-AdG-786580-MACI, SwissMap, and the Einstein Stiftung. 
A.~P. was supported by a fellowship from the Clay Mathematics Foundation.
%D.~Z. was supported by the grant ANR-09-JCJC-0104-01.

This project has received funding from the European Research Council
(ERC) under the European Union's
Horizon 2020 research and innovation program (grant agreement No. 786580).

\section{Curves with an $r$th root} 
\label{Sec:ArtinStack}

\subsection{Artin stacks and Chow cohomology}
\label{ascc}

Let $\ooMM_{g,n}$ denote the smooth Artin stack of prestable curves. 
The Artin stack $\ooMM_{g,n}^Z$ of prestable curves with a line bundle $Z$ of total degree $0$ is obtained from the Picard stack of the universal curve 
$$\pi:\ooCC_{g,n} \to \ooMM_{g,n}\, $$
and is also smooth.
 The Artin stack $\ooMM_{g,n}^Z$
 has a universal curve $$\pi: \ooCC_{g,n}^Z \to \ooMM_{g,n}^Z$$ which 
carries a universal line bundle ${\mathcal Z}$ with sections $s_1, \ldots, s_n$.

Kresch \cite{Kresch} has developed a theory of Chow
 cohomology classes on Artin stacks. 
A basic property is that  given a morphism from a scheme to the stack, 
Kresch's Chow cohomology class on the stack determines a Chow cohomology class on the scheme (comptatible with further pull-backs). 
We describe here a family of Chow cohomology
classes on $\ooMM_{g,n}^Z$.

We first define the set $\G_{g,n}^Z$ of {\em prestable graphs} as follows. 
A prestable graph $\Gamma \in \G_{g,n}^Z$ consists of the data
$$\Gamma=(\V\, ,\ \H\, ,\ \L\, , \ \mathrm{g}:\V \rarr \Z_{\geq 0}\, ,
\ v:\H\rarr \V\, , 
\ \iota : \H\rarr \H,\ \d: \V \to \Z )$$
satisfying the  properties:
\begin{enumerate}
\item[(i)] $\V$ is a vertex set with a genus function $\g:V\to \Z_{\geq 0}$,
\item[(ii)] $\H$ is a half-edge set equipped with a 
vertex assignment $v:\H \to \V$ and an involution $\iota$,
\item[(iii)] $\E$, the edge set, is defined by the
2-cycles of $\iota$ in $\H$ (self-edges at vertices
are permitted),
\item[(iv)] $\L$, the set of legs, is defined by the fixed points of $\iota$ and is
placed in bijective correspondence with a set of $n$ markings,
\item[(v)] the pair $(\V,\E)$ defines a {\em connected} graph
satisfying the genus condition 
$$\nice \sum_{v \in \V} \g(v) + h^1(\Gamma) = g\, ,$$
%\item[(vi)] for each vertex $v$, the stability condition holds:
%$$2\g(v)-2+ \n(v) >0,$$
%where $\n(v)$ is the valence of $\Gamma$ at $v$ including 
%both half-edges and legs,
\item[(vi)] the
degree condition holds:
$$\nice \sum_{v \in \V} \d(v) = 0\, .$$
\end{enumerate}

An automorphism of $\Gamma \in \G^{Z}_{g,n}(X)$ 
consists of automorphisms
of the sets $\V$ and $\H$ which leave invariant the
structures $\L$, $\mathrm{g}$, $v$, $\iota$, and $\d$.
Let $\Aut(\Gamma)$ denote the automorphism group of $\Gamma$.

\begin{definition}
A {\em decorated prestable graph} $[\Gamma,\gamma]$ 
is a prestable graph $\Gamma \in \G_{g,n}^Z$ together with the following
decoration data $\gamma$:
\begin{itemize}
\item each leg $i\in \L$ is decorated with a monomial $\psi_i^a \xi_i^b$.
\item each half-edge $h\in \H\setminus \L$ is decorated
with a monomial $\psi^a_h$,
\item each edge $e\in \E$ is decorated with a monomial $\xi^a_e$,
\item each vertex in $\V$ is decorated with a monomial in the variables 
$\{\eta_{a,b}\}_{a+b \geq 2}$.
\end{itemize}
\end{definition}

%A {\em decorated prestable graph} has, in addition, a monomial $\psi^a \xi^b$ assigned to each leg, a monomial $\psi^a$ assigned to each half-edge that is not a leg, and a monomial in classes $\eta_{a,b}$, $a+b \geq 2$, assigned to each vertex. 

%The labels $\psi$, $\xi$ and $\eta$ have the same meaning as in Notation~\ref{Not:PsiXiEta}. We denote the set of decorated prestable graphs by $\D_{g.n}^S$.

Let $\mathsf{DG}_{g,n}^Z$ be the set of decorated prestable graphs.
Every  $$[\Gamma,\gamma] \in \mathsf{DG}_{g,n}^Z$$ determines a class in the 
Chow cohomology of $\ooMM_{g,n}^Z$ :

\begin{itemize}
\item $\Gamma$ specifies the degeneration of the curve,

\item $\d$ specifies the degree distribution of $Z$,

\item $\psi_i$ corresponds to the cotangent line class,

\item $\xi_i= c_1(s_i^*{\mathcal Z})$,

\item $\xi_e=c_1(s_e^*{\mathcal Z})$ where $s_e$ is the  
node associated to $e$,
\item $\eta_{a,b}=\pi_*\big(c_1(\omega_{\log})^a\, c_1({\mathcal Z})^b\big)$. 
\end{itemize}
We have followed here the 
pattern of Definition \ref{Not:PsiXiEta}.

%We denote the set of decorated prestable graphs by $\D_{g.n}^S$.

More generally, every possibly infinite linear combination of decorated prestable graphs determines  a class in the Chow cohomology of $\ooMM_{g,n}^Z$. 
Indeed, for any morphism  $$B \to \ooMM_{g,n}^Z$$ from a scheme $B$ of
finite type, only a finite number of terms in the linear combination will contribute.
We refer the reader to the Appendix of \cite{GrP2} for the construction of
the product in the Chow cohomology algebra (see also the discussion
in Section \ref{GPGP} below).

\subsection{Twisted curves} \label{Ssec:twistedcurves}

Let $r\geq 1$ be an integer. 
The analog of $\ooMM_{g,n}$ in the context of $r$th roots is
moduli space of  
$\ooMM_{g,n}^r$ of twisted prestable curves constructed in \cite{Ols},
see also \cite{AbrJar,AbrVis}. 
We give a short summary here.

A {\em twisted curve} is a prestable curve with stacky structure at the nodes\footnote{A more complete name is {\em balanced twisted curve}, but we omit the word {\em balanced}, 
since these are the only twisted curves that we consider. 
While stacky structure can also be imposed at the markings of the curve, 
our twisted curves have stacky structure only at the nodes. }. Denote by $\mu_r\subset \C^*$ the group of $r$th roots of unity in the complex plane. The neighborhood of a node in a family of twisted curves is obtained from the family 
$$(x,y) \mapsto z=xy$$ by taking a 
$\mu_r \times \mu_r$ quotient in the source and a $\mu_r$ quotient in the target.

To construct the versal deformation of the node of twisted curve, we start with the
versal deformation
\begin{equation}\label{33992}
  \C^2 \rightarrow  \C\ ,\ \ \ \ (x,y)\mapsto z=xy 
\end{equation}
of the  node of a prestable curve. 
Let $(a,b) \in \mu_r \times \mu_r$ act on $\C^2$ by 
$$(x,y) \mapsto (ax,by)\, ,$$ and let $c \in \mu_r$ act on $\C$ by 
$z \mapsto cz$. These actions commute with \eqref{33992} via the group morphism 
$$\phi: \mu_r \times \mu_r \to \mu_r\, ,\ \ \ \ (a,b) \mapsto c=ab\,. $$
After taking the stack quotient  of both sides of \eqref{33992},
 we obtain a family of twisted curves over $[\C/\mu_r]$ with one stacky $\mu_r$-point at the origin. The fibers of the family over $t \ne 0$ are nonsingular curves isomorphic to $\C^*$. 
The fiber over the origin $t=0$ is the union of the coordinate axes $xy=0$ factored by the kernel of the morphism~$\phi$. 
As soon as a twisted curve acquires a node, the twisted curve
 simultaneously acquires an extra $\mu_r$ group of symmetries (given by the image of $\phi$).

In a family of prestable curves, the neighborhood of a node is 
modeled by the versal deformation 
$$(x,y) \mapsto z = xy\, .$$
Given two line bundles $T_x$ and $T_y$ over a base $B$, 
consider the tensor product  
$$T_z = T_x \otimes T_y\, .$$ 
 We 
construct a family of curves over the total space of $T_z$ over $B$ by 
$$
T_x \oplus_B T_y  \to T_z\, ,\ \ \ 
(b,x,y) \mapsto (b, xy)\, .
$$
Here, $B$ is the boundary divisor and $T_z$ is the normal line bundle to $B$.
We can construct a family of twisted curves over the total space of $T_z$
by applying 
Cadman's $r$th root construction to the zero section $B \subset T_z$. 
In particular, the normal bundle to the locus of nodal twisted
curves is now $(T_x \otimes T_y)^{\otimes (1/r)}$.

A {\em prestable  twisted curve} is a twisted curve with a prestable 
coarsification.
Let $\ooMM^r_{g,n}$ be the moduli space  of prestable twisted curves of genus~$g$ with $n$ marked points. Since
$\ooMM^r_{g,n}$ is obtained from the smooth
Artin stack $\ooMM_{g,n}$ of ordinary prestable curves by applying Cadman's $r$th root construction to the boundary divisor, $\ooMM^r_{g,n}$
is also a smooth Artin stack.
The moduli space $\ooMM^r_{g,n}$  carries {\rm three}  universal curves, see~\cite{Chiodo2}:
\begin{itemize}
\item[(i)] There is the universal twisted prestable curve 
$$\ooCC^r_{g,n} \to \ooMM^r_{g,n}\, .$$ 

\item[(ii)]
There is the fiberwise coarsification 
$\ccoarseC^r_{g,n}$ of $\ooCC^r_{g,n}$. The local model of
$\ccoarseC^r_{g,n}$
is given by 
%the projection of the surface $XY=t^r$ to $t$: it is 
the quotient of the map $$(x,y) \mapsto z= xy$$ by the kernel of the group morphism $\phi$. An  $A_{r-1}$ singularity at the origin is obtained, so
 the universal curve  $\ccoarseC^r_{g,n}$ is singular.
\item[(iii)] There is the
 universal curve $\tC^r_{g,n}$  obtained by resolving the singularities of $\ccoarseC^r_{g,n}$ by a series of blow-ups. 
The resolution of 
the $A_{r-1}$ singularity 
yields a chain of $r-1$ rational exceptional curves over the
origin. The rational curves correspond to the vertices of the $A_{r-1}$ Dynkin diagram, and their intersection points correspond to the edges. 
We call $\tC_{g,n}^r$ the {\em bubbly universal curve}.
\end{itemize}

\subsection{Twisted curves with a line bundle} \label{Ssec:rthroot}

We introduce two more Artin stacks denoted by $\ooMM_{g,n}^{r,Z,\text{triv}}$ and 
$\ooMM_{g,n}^{r,L}$ : 

\vspace{5pt}
\noindent $\bullet\ \ $ The stack $\ooMM_{g,n}^{r,Z,\text{triv}}$ is obtained 
from the stack $\ooMM_{g,n}^{Z}$ of prestable curves with a 
line bundle by applying Cadman's $r$th root construction to the boundary divisors. It is the stack of twisted prestable curves endowed with a 
degree 0 line bundle $Z$ with one extra condition: {\em the stabilizer of every point of the twisted curve acts trivially in the fibers of the line bundle}. 
A Chow cohomology class on $\ooMM_{g,n}^Z$ determines a Chow cohomology class on $\ooMM_{g,n}^{r,Z,{\text{triv}}}$ by pull-back.

\vspace{5pt}
\noindent $\bullet\ $
The stack $\ooMM_{g,n}^{r,L}$ is the stack of twisted prestable curves with a degree 0 line bundle ({\em with no stabilizer conditions}). 
%It can be constructed as follows. First take the Picard stack 
%$$\Pic^0_{\ooCC_{g,n}^r/\ooMM_{g,n}^r}$$ of degree zero line bundles $L$ 
%on twisted curves. 
%This stack comes with a ``dummy'' $\C^*$ action on each line bundle~$L$. The stack $\oMM_{g,n}^{r,L}$ is obtained from $\Pic^0_{\oCC_{g,n}^r/\oMM_{g,n}^r}$ by a rigidification by $\C^*/\mu_r$. In other words, we only keep the dummy action of the finite group $\mu_r$ instead of $\C^*$. In particular, the line bundle $S = L^{\otimes r}$ is completely rigidified.

\vspace{5pt}
These stacks are related by three natural morphisms:
\begin{equation} \label{Eq:p123}
\ooMM_{g,n}^{r,L}\, \stackrel{p_1}{\longrightarrow}\, \ooMM_{g,n}^{r,Z,{\text{triv}}} \, \stackrel{p_2}{\longrightarrow}\, \ooMM_{g,n}^Z\, \stackrel{p_3}{\longrightarrow}\, \ooMM_{g,n}\, .
\end{equation}
The morphism $p_1$ assigns to a pair $(C,L)$ the pair $(C, Z = L^{\otimes r})$, where $C$ is a twisted prestable curve and $L$ a line bundle. The
morphism $p_1$ is \'etale of degree $r^{2g-1}$.
The morphism $p_2$ comes from Cadman's $r$th root construction. The morphism $p_3$ assigns to a pair $(C,S)$ the curve $C$.

While we have taken both $Z$ and $L$ to be of degree 0 in the definitions,
all of the constructions and results of Sections \ref{Sec:ArtinStack}
and \ref{grrrr}  will be valid in case 
$$\deg(Z) = r \deg(L)\, .$$

\subsection{Commutative diagram}
Let $X$ be a nonsingular projective variety with line bundle
$S\rightarrow X$. Let  $A=(a_1,\ldots, a_n)$ be a vector
of integers which satisfy
$$\sum_{i=1}^n a_i = \int_\beta c_1(S)$$
for $\beta\in H_2(X,\mathbb{Z})$.

The moduli 
space $\oM_{g,A,\beta}^r(X,S)$ of stable maps
$$f:(C,x_1,\ldots,x_n) \to X$$
endowed with an $r$th root of the degree 0 line bundle 
$$f^*S\Big(-\sum_{i=1}^n a_i x_i\Big)\, $$
plays a central role in the proof of Theorem \ref{Thm:main}.
 The moduli space 
$\oM_{g,A,\beta}^r(X,S)$
is defined as the  fiber product of the following two maps:

\vspace{5pt}
\noindent (i) \ $\pi_Z: \oM_{g,n,\beta}(X)\to \ooMM_{g,n}^Z$
assigns to a stable map $f: C \to X$ the pair 
$$
\Big(C, f^*S\Big(-\sum_{i=1}^n a_i x_i\Big)\Big)\, .
$$

\vspace{5pt}
\noindent (ii)\
$\epsilon : \ooMM_{g,n}^{r,L} \to \ooMM_{g,n}^Z$ is the composition 
$\epsilon = p_2 \circ p_1$.

\vspace{5pt}
\noindent The moduli space $\oM_{g,A,\beta}^r(X,S)$ is the fiber product of 
$\pi_Z$ and $\epsilon$:

\vspace{-3.0em}
\begin{equation} \label{Eq:CommDiag}
\setlength{\unitlength}{1em}
%\begin{center}
\begin{picture}(13,7)(0,3)
\put(0,6){$\oM_{g,A,\beta}^r(X,S)$}
\put(2,0){$\ooMM_{g,n}^{r,L}$}
\put(9,6){$\oM_{g,n,\beta}(X)$}
\put(10.2,0){$\ooMM_{g,n}^Z\ ,$}
\put(5,.35){\vector(1,0){4.3}}
\put(6.5,6.35){\vector(1,0){2}}
\put(3,5.1){\vector(0,-1){3.5}}
\put(10.9,5.1){\vector(0,-1){3.5}}
\put(6.9,0.6){\mbox{\small $ \epsilon$ }}
\put(7.3,6.5){\mbox{\small $\epsilon$ }}
\put(9.7,3.3){\mbox{\small $\pi_Z$ }}
\put(1.8,3.3){\mbox{\small $\pi_L$ }}
\end{picture}
%\end{center}
\end{equation}

%\vspace{4\unitlength}
\vspace{35pt}

\noindent
where we denote top arrow also by $\epsilon$.

%both forgetful maps that forget the line bundle~$L$. 
%Both are finite of degree~$r^{2g-1}$. The map 

\subsection{The pull-back map $\pi_Z^*$} \label{Ssec:pullback}

We describe here the pull-back map $\pi_Z^*$ 
for Chow cohomology classes defined by decorated graphs.
Let $$[\Gamma,\gamma] \in \D\G_{g,n}^Z$$ be a decorated prestable graph
representing a Chow cohomology class in $\ooMM_{g,n}^Z$
following the conventions of Section \ref{ascc}.
% while $(\Gamma, w)$ represents one in~$\oMM_{g,n}^{r,L}$.

\begin{lemma} \label{Lem:pullback1}
The pull-back $\pi_Z^*[\Gamma,\gamma]$  
%$\pi_L^*[\Gamma,w]$ 
is obtained in terms of 
Chow cohomology classes of decorated $X$-valued stable 
graphs by applying the following procedure:
\begin{itemize}
    \item Replace the degree $\d(v) \in \Z$ of each vertex 
with effective classes $$\beta(v) \in H_2(X,\Z)$$ satisfying 
$$\int_{\beta(v)} c_1(S) -\sum_{i \vdash v} a_i = \d(v)\,  ,$$ 
where the sum is over the legs $i$ incident to $v$. 
Sum over all choices of~$\beta(v)$.
 \item Replace each $\xi_i$ with $\xi_i+a_i\psi_i$.
    \item  Replace each class $\eta_{0,b}$
at each vertex~$v$ with 
    $$
    \eta_{0,b}-\sum_{i \vdash v} \sum_{k=1}^b
    \binom{b}{k} a_i^k \psi_i^{k-1} \xi_i^{b-k},
    $$
    where the first sum is again over the legs $i$ incident to~$v$.
   
\end{itemize}
 All other decorations are kept the same.
\end{lemma}

\paragraph{Proof.} Given a stable map $f:C \to X$, the degree of 
$f^*S(-\sum_{i=1}^n a_i x_i)$ on the component of the curve 
$C$ corresponding a vertex $v$ of the dual graph equals  
$$\int_{\beta(v)} c_1(S)-\sum_{i \vdash v} a_i\, ,$$
 which justifies the first operation. 

Recall the divisor $D_i$ corresponds to the $i$th section, 
$$ D_i \subset \ooCC_{g,n}^Z \stackrel{\pi}{\rightarrow} \ooMM_{g,n}^Z\, ,$$
and $\omega_\log = \omega_\pi(\sum_{i=1}^n D_i)$. 
 The first Chern class of $f^*S(-\sum_{i=1}^n a_i x_i)$ on the universal curve equals 
$$c_1(S) - \sum_{i=1}^n a_i D_i\, .$$ 
The pull-back 
%DZ: \pi_S --> \pi_Z
$\pi_Z^*(\eta_{a,b})$ %or %$\pi_L^* \eta_{a,b}$ 
is the push-forward from the universal curve of the product
\begin{equation}\label{ff566}
K^a \Big(c_1(S) - \sum_{i=1}^n a_i D_i\Big)^b\, ,
\end{equation}
where
$K = c_1(\omega_\log)$.
Since $K D_i =0$, the product \eqref{ff566} is equal to 
$K^a \xi^b$ if $a>0$. Hence, for $a>0$, 
%DZ: \pi_S --> \pi_Z
$$\pi_Z^*(\eta_{a,b}) = \eta_{a,b}\, ,$$
% $\pi_L^*\eta_{a,b} = \eta_{a,b}$ for $a>0$.
 In  case $a=0$, we expand $(c_1(S) - \sum_{i=1}^n a_i D_i)^b$ and take 
the push-forward from the universal curve to the moduli space.
We find
$$
    \eta_{0,b}-\sum_{i \vdash v} \sum_{k=1}^b
    \binom{b}{k} a_i^k \psi_i^{k-1} \xi^{b-k}_i
$$
in the notation of decorated $X$-valued stable graphs.
\qed

\subsection{Chow cohomology classes on 
 $\oM_{g,A,\beta}^r(X,S)$ and $\ooMM_{g,n}^{r,L}$}

\subsubsection{Overview}
Recall the commutative diagram~\eqref{Eq:CommDiag}:

\vspace{-40pt}
\begin{equation*} 
\setlength{\unitlength}{1em}
%\begin{center}
\begin{picture}(13,7)(0,3)
\put(0,6){$\oM_{g,A,\beta}^r(X,S)$}
\put(2,0){$\ooMM_{g,n}^{r,L}$}
\put(9,6){$\oM_{g,n,\beta}(X)$}
\put(10.2,0){$\ooMM_{g,n}^Z$}
\put(5,.35){\vector(1,0){4.3}}
\put(6.5,6.35){\vector(1,0){2}}
\put(3,5.1){\vector(0,-1){3.5}}
\put(10.9,5.1){\vector(0,-1){3.5}}
\put(6.9,0.6){\mbox{\small $ \epsilon$ }}
\put(7.3,6.5){\mbox{\small $\epsilon$ }}
\put(9.7,3.3){\mbox{\small $\pi_Z$ }}
\put(1.8,3.3){\mbox{\small $\pi_L$ }}
\end{picture}
%\end{center}
\end{equation*}

%\vspace{4\unitlength}
\vspace{40pt}

\noindent We now define Chow cohomology classes via decorated graphs on
$$\oM_{g,A,\beta}^r(X,S) \ \  \text{and} \ \ \ \ooMM_{g,n}^{r,L}$$
and describe the pull-back map $\pi^*_L$.
Except for the additional data recording the twisted structure,
the discussion is almost identical to our treatment of
$$\pi_Z: \oM_{g,n,\beta}(X,S) \rightarrow \ooMM_{g,n}^Z\, .$$

\subsubsection{The moduli space $\oM_{g,A,\beta}^r(X,S)$} \label{www111}
We define the set $\G^r_{g,A,\beta}(X,S)$ of {\em $X$-valued $r$-twisted
  stable graphs} as
follows. A graph $\Gamma \in \G^r_{g,A,\beta}(X,S)$ consists of the data
$$\Gamma=(\V\, ,\ \H\, ,\ \L\, , \ \mathrm{g}, 
\ v\, , 
\ \iota\, ,\ \beta: \V \to H_2(X,\Z), \tw:\H \to \{0,\ldots,r-1\})$$
satisfying the properties:
\begin{enumerate}
\item[(i-vii)] exactly as for
   $X$-valued 
  stable graphs in Section \ref{xvalsg},
\item[(viii)] the twist conditions hold:
\begin{itemize}
\item[($\L$)]  $\forall i\in \L\ \Longrightarrow \ \tw(i)=0\, ,$
\item[($\E$)] $\forall e = (h', h'')\in \E \ \Longrightarrow \  
 \tw(h') + \tw(h'') = 0 \bmod r\, ,$
\item[($\V$)]
$\forall v\in \V \ \Longrightarrow \ 
\sum\limits_{\nu(h)=v} \tw(h) = \int_{\beta(v)} c_1(S) -\sum_{i \vdash v} a_i\ \bmod r\, .$
\end{itemize}  
\end{enumerate}
The line bundle $S\rightarrow X$ and vector $A$
 only appear in property (viii).

The universal curve $$\pi:\ooC^r_{g,A,\beta} \to \oM^r_{g,A,\beta}(X,S)$$ carries the log relative dualizing sheaf $\omega_{\log}$ 
and the pull-back $f^*S$  of the line bundle $S$ via the
universal map,
$$f: \ooC^r_{g,A,\beta} \rightarrow X\, .$$
%Let  $K = c_1 (\omega_\log)$ and
%$\xi = c_1(f^*(S))$.
Following Section \ref{Ssec:PsiXiEta}, we
define Chow cohomology classes
on $\oM^r_{g,A,\beta}(X,S)$:
\begin{itemize}
    \item $\psi_i = c_1(s_i^* \omega_\pi)\, $,
    \item $\xi_i = c_1(s_i^* f^*S)\, $,
    \item $\eta_{a,b} = \pi_*\left(c_1(\omega_\log)^a\, c_1(f^*S)^b\right)\, $.
\end{itemize}

\begin{definition}
A {\em decorated $X$-valued $r$-twisted stable graph} $[\Gamma,\gamma]$ 
%of type $(g,n,\beta)$ 
is an 
$X$-valued stable graph $\Gamma \in \G^r_{g,A,\beta}(X)$ together with the following
decoration data $\gamma$:
\begin{itemize}
\item each leg $i\in \L$ is decorated with a monomial $\psi^a_i \xi_i^b$,
\item each half-edge $h\in \H\setminus \L$ is 
decorated
with a monomial $\psi^a_h$,
\item each edge $e\in \E$ is decorated with a monomial $\xi^a_e$,
\item each vertex in $\V$ is decorated with a monomial in the variables 
$\{\eta_{a,b}\}_{a+b \geq 2}$.
\end{itemize}

\end{definition}

Let $\D\G^r_{g,A,\beta}(X,S)$ be the set of decorated $X$-valued $r$-twisted
stable graphs.
Every  $$[\Gamma,\gamma] \in \mathsf{DG}_{g,A,\beta}^{r}(X,S)$$ determines a class in the
Chow cohomology of $\oM_{g,A,\beta}^{r}(X,S)$ :

\begin{itemize}
\item $\Gamma$ specifies the degeneration of the curve with twisted
structure $w$ at the nodes on the $r$th root of 
$$f^*S\Big(-\sum_{i=1}^n a_i x_i\Big)\, ,$$

\item $\beta$ specifies the curve class distribution of the map $f$,

\item the decorations $\psi_i$, $\xi_i$, $\xi_e$, and $\eta_{a,b}$
specify Chow cohomology classes, 
$$\xi_e = c_1(s_e^*f^*S)\, .$$
\end{itemize}
%We have followed the 
%notation of Definition \ref{Not:PsiXiEta}.

\subsubsection{The Artin stack $\ooMM_{g,n}^{r,L}$} \label{www222}
Let $\G_{g,n}^{r,L}$ be the set of {\em $r$-twisted prestable graphs}
 defined as follows.
A  graph $$\Gamma \in \G_{g,n}^{r,L}$$
 consists of the data
$$\Gamma=(\V\, ,\ \H\, ,\ \L\, , \ \mathrm{g}\, ,
\ v \, , 
\ \iota \, , \ \d: \V \to \Z\, , \tw:\H\rightarrow \{0,\ldots,r-1\} )$$
satisfying the properties:
\begin{enumerate}
\item[(i-vi)]
exactly as for
   prestable
  graphs in Section \ref{ascc},  
\item[(vii)] the twist conditions hold:
\begin{itemize}
\item[($\L$)] $\forall i\in \L\ \Longrightarrow  \ \tw(i)=0\, ,$
\item[($\E$)] $\forall e = (h', h'')\in \E \ \Longrightarrow \  
 \tw(h') + \tw(h'') = 0 \bmod r\, ,$
\item[($\V$)]
$\forall v\in \V \ \Longrightarrow \ 
\sum\limits_{\nu(h)=v} \tw(h) = \d(v) \bmod r\, .$
\end{itemize}
\end{enumerate}

\begin{definition}
A {\em decorated $r$-twisted prestable graph} $[\Gamma,\gamma]$ 
is a graph $$\Gamma \in \G_{g,n}^{r,L}$$ together with the following
decoration data $\gamma$:
\begin{itemize}
\item each leg $i\in \L$ is decorated with a monomial $\psi^a_i \xi^b_i$,
\item each half-edge $h\in \H\setminus \L$ is decorated
with a monomial $\psi^a_h$,
\item each edge $e\in \E$ is decorated with a monomial $\xi^a_e$,
\item each vertex in $\V$ is decorated with a monomial in the variables 
$\{\eta_{a,b}\}_{a+b \geq 2}. $
\end{itemize}
\end{definition}

%A {\em decorated prestable graph} has, in addition, a monomial $\psi^a \xi^b$ assigned to each leg, a monomial $\psi^a$ assigned to each half-edge that is not a leg, and a monomial in classes $\eta_{a,b}$, $a+b \geq 2$, assigned to each vertex. 

%The labels $\psi$, $\xi$ and $\eta$ have the same meaning as in Notation~\ref{Not:PsiXiEta}. We denote the set of decorated prestable graphs by $\D_{g.n}^S$.

Let $\mathsf{DG}_{g,n}^{r,L}$ be the set of decorated $r$-weighted
prestable graphs.
Every  $$[\Gamma,\gamma] \in \mathsf{DG}_{g,n}^{r,L}$$ determines a class in the
Chow cohomology of $\ooMM_{g,n}^{r,L}$ :

\begin{itemize}
\item $\Gamma$ specifies the degeneration of the curve with twisted
structure $w$ at the nodes on $L$,

\item $\d$ specifies the degree distribution of $L^{\otimes r}$,

\item $\psi_i$ corresponds to the cotangent line class,

\item $\xi_i= c_1(s_i^*L^{\otimes r})$,

\item $\xi_e=c_1(s_e^*L^{\otimes r})$ where $s_e$ is the  
node associated to $e$,

\item $\eta_{a,b}=\pi_*\big(c_1(\omega_{\log})^a\, c_1(L^{\otimes r})^b\big)$.
\end{itemize}

\subsubsection{The pull-back map $\pi^*_L$}
Let $[\Gamma,\gamma] \in \D\G_{g,n}^{r,L}$ be a decorated prestable graph
representing a Chow cohomology class in $\ooMM_{g,n}^{r,L}$.
The pull-back
$\pi^*_L$ along
$$\pi_L: 
\oM_{g,A,\beta}^r(X,S) \to \ooMM_{g,n}^{r,L}$$
is computed by exactly the same rules governing $\pi^*_Z$.
The proof is identical.

\begin{lemma} \label{Lem:pullback}
The pull-back $\pi_L^*[\Gamma,\gamma]$   
is obtained in terms of 
Chow cohomology classes of decorated $X$-valued $r$-twisted stable 
graphs by applying the following procedure:
\begin{itemize}
    \item Replace the degree $\d(v) \in \Z$ of each vertex 
with all effective classes $$\beta(v) \in H_2(X,\Z)$$
 satisfying 
$$\int_{\beta(v)} c_1(S) -\sum_{i \vdash v} a_i = \d(v)\,  ,$$ 
where the sum is over the legs $i$ incident to $v$. 
Sum over all choices of~$\beta(v)$.
\item Replace each $\xi_i$ with $\xi_i+a_i\psi_i$.
    \item  Replace each class $\eta_{0,b}$
at each vertex~$v$ with 
    $$
    \eta_{0,b}-\sum_{i \vdash v} \sum_{k=1}^b
    \binom{b}{k} a_i^k \psi_i^{k-1} \xi_i^{b-k},
    $$
    where the sum is again over the legs $i$ incident to~$v$.
\end{itemize}
All other decorations are kept the same.
\end{lemma}

\subsection{Multiplication in the Chow cohomology of $\ooMM_{g,n}^{r,L}$}
\label{GPGP}
The product in Chow cohomology of
the classes of two decorated $r$-twisted prestable graphs
in $\D\G_{g,n}^{r,L}$ 
is defined in a very similar way to the product of stable graphs carefully described in the appendix of~\cite{GrP2}. We briefly sketch the construction and
highlight the differences.

An {\em edge contraction} in an $r$-twisted prestable graph is defined
in the natural way: 
\begin{itemize}
\item if the edge contraction merges two vertices, the
corresponding genera and degrees are summed,
\item  if the contracted  edge is a loop, the genus of the base vertex increases by~1 and the degree remains unchanged.
\end{itemize}
The total degree and twisting conditions are still satisfied after 
edge contraction in  an $r$-twisted prestable graph.
%The decorations are also canonically merged in an edge contraction.

We define the product of two decorated $r$-twisted prestable graphs 
$$[\Gamma_A,\gamma_A]\, ,\ [\Gamma_B,\gamma_B]\, \in \D\G_{g,n}^{r,L}\, $$
as follows. The product is a (possibly infinite) linear combination  
of decorated $r$-twisted prestable graphs.

We first consider prestable graphs
$$\Gamma\in \G_{g,n}^{r,L}$$
 with edges colored by $A$, $B$, or both $A$ and $B$, 
satisfying the conditions:
\begin{itemize}
\item[(i)]
after contracting the edges not colored $A$, we
 obtain $\Gamma_A$, 
\item[(ii)]
after contracting the edges not colored $B$, we obtain $\Gamma_B$.
\end{itemize}
For each such $\Gamma$, we add decorations by the following rules.
\begin{itemize}
\item The monomials $\psi^a \xi^b$ on the legs of the graph $\Gamma$ are obtained by multiplying the corresponding leg  monomials on the graphs $\Gamma_A$ and $\Gamma_B$. 
\item The monomial $\psi^a$ on a half-edge colored $A$ only or colored $B$ only 
is inherited from the graph $\Gamma_A$ or $\Gamma_B$ respectively. 
On an edge $e=(h',h'')$ colored both $A$ and $B$, we
take the product of the monomials on the corresponding edge in the graphs $\Gamma_A$ and $\Gamma_B$  and include an
 extra factor $$-\frac{1}{r}(\psi_{h'} + \psi_{h''})$$ 
corresponding to the excess intersection. 

\item
%The factors of $\xi$ on the edges are obtained by multiplying
%the corresponding edge factors of $\Gamma_A$ and $\Gamma_B$ respectively.
  The monomial $\xi^a$ on an edge colored $A$ only or colored $B$ only is inherited from the graph $\Gamma_A$ or $\Gamma_B$ respectively. On an edge colored both
  $A$ and $B$, we take the product of the monomials on the corresponding edge in the graphs $\Gamma_A$ and $\Gamma_B$.

\item The factors $\eta_{a,b}$ of the monomials assigned to each vertex $v$ of $\Gamma_A$ (and $\Gamma_B$) are 
%DZ: ``split'' replaced with ``distributed'', a clarifying sentence added.
distributed in all possible ways among the vertices 
which collapse to $v$ as $\Gamma$ is contracted to $\Gamma_A$ (and $\Gamma_B$). In other words, each factor is assigned to a unique vertex, and we sum over all such assignments. Each vertex of $\Gamma$ is then marked with two monomials in the variables $\eta_{a,b}$, which we multiply together.
\end{itemize}
The product 
$[\Gamma_A,\gamma_A]\cdot [\Gamma_B,\gamma_B]$
is then the sum over all decorated $r$-twisted prestable
graphs $[\Gamma,\gamma]$ produced by the above construction.

%It is important to note that the product of two infinite linear combinations of decorated $r$-twisted prestable graphs is well defined since the coefficient of each graph in the product only involves a finite number of graphs in the factors. The product rule transforms the vector space spanned by the decorated 
%$r$-twisted prestable graphs into an algebra which
%agrees with the product in the Chow cohomology of

It is important to note that 
a product of two decorated $r$-twisted prestable graphs can be an
{\em infinite} linear combination of decorated $r$-twisted prestable graphs.
For instance, if we take the square of
the graph with a single vertex of degree~0 and a single loop (and no decorations),
we obtain, among other terms, a sum over all graphs with two vertices of degrees
$d$ and $-d$ connected by two edges. Since the integer $d$ can be chosen arbitrarily,
the result is an infinite linear combination.
However, the product of two  decorated $r$-twisted prestable graphs
or even of infinite linear combinations of such graphs,
is always well-defined because the coefficient of each graph in the product
only involves a finite number of graphs in the factors. Therefore,
the product rule transforms the vector space of possibly infinite linear combinations of decorated $r$-twisted prestable graphs into an algebra which agrees with the product in the Chow cohomology of $\ooMM_{g,n}^{r,L}$.

\section{GRR for the universal line bundle} \label{grrrr}

\subsection{Universal  bundles on universal curves over $\ooMM_{g,n}^{r,L}$}

\label{pp234}
The universal twisted curve 
$$\pi:\ooCC_{g,n}^{r,L} \to \ooMM_{g,n}^{r,L}$$ carries a universal 
line bundle $\cL$. 
Consider a node in a singular fiber of the universal curve. 
The kernel of the group morphism $\phi$ of Section \ref{Ssec:twistedcurves}
acts on the fiber of~$\cL$ over the node. 
If the generator $$(1,-1) \in \mu_r \times \mu_r$$ of the kernel acts on the fiber of 
$L$ at the node by $e^{2 \pi a /r}$, 
we assign the remainders $a$ and $-a$ mod~$r$ to the branches of the curve at the node.
Every node in the universal curve therefore acquires a {\em type}\/: 
a pair of remainders mod~$r$ assigned to the branches meeting at the node (with vanishing sum mod~$r$). 
For the line bundle $$\ZZ=\cL^{\otimes r}\,,$$
 the action of the kernel of $\phi$ is always trivial. 

The push-forward $\coarseL$ of the sheaf of invariant sections of $\cL$ to the coarse universal curve 
$\mathsf{C}_{g,n}^{r,L}$, is a rank~1 torsion free sheaf which is described in detail in \cite{ChiZvo}. 
The  push-forward of the sheaf of sections of~$\cL^{\otimes r}$  to 
$\mathsf{C}_{g,n}^{r,L}$
is just a line bundle $\widetilde{\ZZ}$.

On the bubbly universal curve  $\tC_{g,n}^{r,L}$,
 we may pull-back the line bundle $\widetilde{\ZZ}$  from the coarse  curve~
$\mathsf{C}_{g,n}^{r,L}$. The situation is more interesting for $\cL$. 
Chiodo~\cite{Chiodo2} has proven that there exists a line bundle
$$\tL\to \tC_{g,n}^{r,L}$$ for which
 the push-forward of the associated sheaf of sections to the coarse  curve 
$\mathsf{C}_{g,n}^{r,L}$ is~$\coarseL$. However, instead of the simple isomorphism 
$$\ZZ = \cL^{\otimes r} \ \ \text{on}\ \ \ooCC_{g,n}^{r,L}\,,$$ we have the more
complicated relation
$$
\widetilde{\ZZ} = \tL^{\otimes r}(D) \ \ \text{on}\ \ \tC_{g,n}^{r,L}
 \, ,
$$
where $D$ is a linear combination of the exceptional divisors 
of the desingularization 
$$\tC_{g,n}^{r,L} \to \mathsf{C}_{g,n}^{r,L}\, .$$
 More precisely, a node of type $$(a,b) \ \ \text{with} \ \ a+b = 0 \pmod r\, ,$$
 gives rise to a chain of $r-1$ rational curves in the fiber of the bubbly universal curve. 
These rational curves appear in~$D$ with coefficients 
$$
a, \; 2a,\;  \dots, \; (b-1)a, \; ab, \; (a-1)b, \; \dots, \; 2b, \; b,
$$
see~\cite{Chiodo2}. 

\subsection{Polynomial classes}

%DZ: new paragraph.
Ehrhart's theory states that if we take an $n$-dimensional polytope $\Delta$ with integer vertices and a polynomial $P(x_1, \dots, x_n)$ in $n$ variables, then the sum of values of $P$ over the integer points inside $r \Delta$ is a polynomial in~$r$. Actually, the same claim holds for a polynomial $P(x_1, \dots, x_n;r)$ in $n+1$ variables depending explicitly on~$r$. Similarly to the discussion in the Appendix of \cite{JPPZ}, we will use Ehrhart's theory to prove that a family of cohomology classes on $\ooMM_{g,n}^{r,L}$ projected to $\ooMM_{g,n}^Z$ forms a cohomology-valued Laurent polynomial in~$r$.

Consider a family of polynomials
$$\Big\{\, P_\Gamma \in \C[\tw_h, r, r^{-1}]\, \Big\}_{\Gamma\in \D\G^Z_{g,n}}\, .$$
For each graph $  \Gamma\in \D\G^Z_{g,n}$,
$P_\Gamma$ is a polynomial in the variables $\tw_h$, where $h$ runs over the half-edges of~$\Gamma$ which {\em are not legs}
(and $P_\Gamma$ is a
 Laurent polynomial in~$r$). 
%DZ: added sentence.
The formal variables $\tw_h$ play the role of the variables $x_i$ in the previous paragraph.
 
The family $P_\Gamma$  determines a family of Chow cohomology classes on 
$\ooMM_{g,n}^{r,L}$ for all~$r$:
$$
\alpha = \sum_{\Gamma} \sum_{\tw} P_\Gamma\left(\tw(h), r,r^{-1}\right) \cdot [\Gamma, \tw]\, ,
$$
where the summation is over all $r$-twistings $\tw$ of
all decorated prestable graphs $\Gamma\in \D\G_{g,n}^{Z}$ 
which, equivalently, is the set $\D\G_{g,n}^{r,L}$. 
%DZ: added sentence.
Note that we have substituted
the value of the twist $\tw(h) \in \{0, \dots, r-1 \}$ in place of the formal variable $\tw_h$.

We call such families of classes {\em polynomial}. We let
$$
\val(\alpha) = \inf_{\Gamma \in \D\G_{g,n}^Z} \Big[\val_r(P_\Gamma) - |\E(\Gamma)|
\Big]\, .
$$

\begin{proposition} \label{Prop:product}
  If $\alpha$ and $\beta$ are two polynomial families of classes of valuations $\val(\alpha)$ and $\val(\beta)$, then the product\,{\footnote{The
      product of two families
      of Chow cohomology classes is defined
      by taking the products of the corresponding Chow cohomology
      classes on 
$\ooMM_{g,n}^{r,L}$ for all~$r$.}}
      %DZ: removed ``of''.
       $\alpha \beta$ is a polynomial family of classes satisfying
$$\val(\alpha\beta) \geq \val(\alpha) + \val(\beta)\, .$$
\end{proposition}

\paragraph{Proof.}
Let $\Gamma$ be a decorated prestable graph and $\tw$ an $r$-twisting.
 The coefficient of $\Gamma$ in the product $\alpha \beta$ is a finite linear combination of coefficients of contractions of $\Gamma$ in $\alpha$ and $\beta$.
 More precisely, the coefficient is a sum over all ways to label the 
edges of $\Gamma$ with $A$, $B$, or $AB$ and  
to split every monomial on a leg, edge, or vertex  of $\Gamma$ into 
a product of two factors labeled $A$ and $B$ while keeping a factor 
$$-\frac{1}{r}(\psi'+\psi'')$$ aside for every edge labeled $AB$. For every such labelling,
the contribution to the coefficient of $\Gamma$ is the product of the
two corresponding polynomials in $\tw(h)$, $r$, and $r^{-1}$
 times a factor of $r^{-1}$ for each $AB$-edge. 
Therefore the coefficient of $\Gamma$ is also a polynomial in  
$\tw(h)$, $r$, and $r^{-1}$. 

The lowest possible order in~$r$ is given by summing the following three degrees:
\begin{itemize}
\item $\val(\alpha) + (\mbox{number of } A\mbox{-edges} \ \text{and} \ AB\mbox{-edges})$
  from the coefficient of the $A$-contraction of $\Gamma$ in $\alpha$,
\item $\val(\beta) + (\mbox{number of } B \mbox{-edges}
\ \text{and} \ AB\mbox{-edges})$
  from the coefficient of the $B$-contraction of $\Gamma$ in $\beta$,
    \item $- (\mbox{number of } AB \mbox{-edges})$ from the excess intersection factor $$-\frac{1}{r}(\psi' + \psi'')\, .$$
\end{itemize}
The sum of these three degrees is, indeed, $\val(\alpha) + \val(\beta)+ |E(\Gamma)|$, so the valuation of $\alpha \beta$ is the sum of valuations of $\alpha$ and $\beta$ (or larger if the lowest degree terms cancel out). \qed

\begin{proposition} \label{Prop:pushforward}
Let $\alpha$ be a polynomial family of classes in the strata algebra of $\ooMM_{g,n}^{r,L}$. Let $\beta$ be the push-forward of $\alpha$ to the strata algebra of $\ooMM_{g,n}^Z$. Then, for any decorated prestable graph $\Gamma \in \D\G_{g,n}^Z$, the coefficient of $\Gamma$ in $\beta$ is a Laurent polynomial in~$r$ of valuation at least $\val(\alpha)+2g-1$,
 for all  sufficiently large $r$.
\end{proposition}

\paragraph{Proof.} 
%DZ: added sentence.
The projection from $\ooMM_{g,n}^{r,L}$ to $\ooMM_{g,n}^Z$ sums over all possible twists, which is the analog of summing the values of a polynomial over integer points in a polytope in Ehrhart's theory.

For every $r$-twisting $\tw$ of $\Gamma$,
 the push-forward from the stratum $[\Gamma, \tw]$ of $\ooMM_{g,n}^{r,L}$ to the stratum $\Gamma$ of $\ooMM_{g,n}^Z$ has degree
$$
r^{\sum_{v \in V(\Gamma)} (2 \g(v)-1)}
= r^{2g - 2 h^1(\Gamma) - |V(\Gamma)|}\, 
$$
see the proof of \cite[Corollary 4]{JPPZ}.

The coefficient of $\Gamma$ in $\beta$ is obtained by 
summing the coefficients of $[\Gamma, \tw]$ in $\alpha$ over all 
$r$-twistings $\tw$ of $\Gamma$. According to Proposition $3''$ of the Appendix of~\cite{JPPZ}, the sum is a Laurent polynomial in~$r$, for $r$ 
sufficiently large, with valuation  at least 
$$
\val(\alpha) + |E(\Gamma)| + h^1(\Gamma)\, .$$
Multiplying the Laurent polynomial by the degree of the push-forward map,
 we obtain, again, a Laurent polynomial in~$r$ with valuation at least
$$
\val(\alpha) + |E(\Gamma)| + h^1(\Gamma) + 2g - 2 h^1(\Gamma) - V(\Gamma) = \val(\alpha) + 2g-1\, 
$$
as claimed.
\qed

\subsection{Applying the GRR formula to $\tC_{g,n}^{r,L} \to \ooMM_{g,n}^{r,L}$}
\label{Ssec:GRR}

The bubbly universal curve
$$\pi: \tC_{g,n}^{r,L} \to \ooMM_{g,n}^{r,L}$$ is  representable and
proper over $\ooMM_{g,n}^{r,L}$. We may  apply the Grothendieck-Riemann-Roch (GRR) formula
to compute $R\pi_* \tL$,
$$
\ch \left(R\pi_* \tL\right) = \pi_*\left(\ch(\tL)\cdot \td(\pi)\right).
$$

As before, let $D_i$ be the class of the
divisor of $\tC_{g,n}^{r,L}$ corresponding to the $i$th section.
Let
$$K = c_1(\omega_\log) = c_1\Big(\omega\Big(\sum_{i=1}^n D_i\Big)\Big)\, , \ \ \ \xi = c_1\big(\widetilde{\ZZ}\big)\,.
$$
Following the notation of Section \ref{pp234}, let $D$ be defined by 
$$
\widetilde{\ZZ} = \tL^{\otimes r}(D)\, .
$$

Let $j$ be the double covering of the locus of nodes in the singular fibers,
$$j: \Delta \to \tC_{g,n}^{r,L}\, .$$
The sheets of the covering correspond to the two ways of numbering the branches of the curve at the node. The domain $\Delta$ carries two cotangent
line bundles
corresponding to the two branches of the node.
We denote by $\nu_1$ and $\nu_2$ first Chern classes of the
two cotangent lines.

Using the above notation, we have:
$$
\ch(\tL) =  e^{\xi/r}e^{-D/r},
$$
$$
\td(\pi) = \frac{K}{e^K-1} 
\prod_{i=1}^n \frac{D_i}{1-e^{-D_i}} 
\left( 1 + j_*\,\frac{1}{2} \sum_{m \geq 1} \frac{B_{2m}}{(2m)!} \frac{\nu_1^{2m-1} + \nu_2^{2m-1}}{\nu_1 + \nu_2}
\right),
$$
see~\cite{Chiodo2,FP-H}.

\begin{lemma} \label{Lem:GRR}
We have 
$$
\pi_*\left(\ch(\tL) \td(\pi)\right) = 
$$
$$
\pi_* \left[
e^{\xi/r}\frac{K}{e^K-1} 
\prod_{i=1}^n \frac{D_i}{1-e^{-D_i}}
\; + \;
e^{\xi/r-D/r} \left( 1 + j_*\, \frac{1}{2}\sum_{m \geq 1} \frac{B_{2m}}{(2m)!} \frac{\nu_1^{2m-1} + \nu_2^{2m-1}}{\nu_1 + \nu_2}
\right)
\right].
$$
\end{lemma}

\paragraph{Proof.}
The following intersections vanish in
$ \tC_{g,n}^{r,L}$
by the definition of $\omega_{\log}$ and
the fact that the markings are disjoint from the nodes (before the bubbly resolution
and, therefore, also after the bubbly resolution):
$$K\cdot \Delta = D_i\cdot \Delta=D_i\cdot D=0\, . $$
Moreover, since the the bubble curve is crepant (see \cite{Chiodo2}), we have
$$ K \cdot D=0\, .$$
The push-forward of the constant term vanishes $$\pi_*(1)=0\, .$$
% from the first term with any class $D$ or $\Delta$ from the second term vanishes. Indeed, $D \cdot K$  and $\Delta \times K$ vanish, because the singularity resolution in the bubbly curve is crepant, see~\cite{Chiodo2}; $D \cdot D_i$  and $\Delta \cdot D_i$ vanish  because their geometric intersection is empty; $D \times \xi$ and $\Delta \cdot \xi$ vanish, because the line bundle $S$ is a pull-back from the coarse curve, so that its restriction to $D$ and $\Delta$ is trivial.
The formula is then easily obtained from the above vanishing from GRR.
\qed

\subsection{GRR and polynomiality in~$r$}
\label{Ssec:GRRpoly}

\begin{proposition} \label{Prop:character}
The GRR formula for the Chern character $\ch_k(R\pi_* \tL)$ is a polynomial family with valuation $-(k+1)$. Its coefficients of lowest degree in~$r$ equal
$$
r^{-(k+1)}\, \eta_{0,k+1} \; + \; \; r^{-k} \!\!\!\!\! \sum_{\Gamma \;
\mbox{\tiny \rm with one edge } (h,h')
} \
%\frac{w(h)^{k+1}}{(k+1)!} \sum_{i+j= k-1} (-1)^i \psi_h^i \psi_{h'}^j.
\sum_{i+j+m =k-1}
\frac  { \xi_e^m }  { m! }
\frac  {  \tw(h)^{k+1-m}  }  {  (k+1-m)!  }
(-1)^i  \psi_h^i  \psi_{h'}^j\, .
$$
\end{proposition}

\paragraph{Proof.} The first term in Lemma~\ref{Lem:GRR} only involves decorated graphs with a single vertex  and with a single
 (trivial) weighting. The dependence on $r$ in the first term
is only through $e^{\xi/r}$. The degree~$k$ part of the class, obtained by 
$\pi$ push-forward of the degree~$k+1$ part, is a Laurent polynomial in~$r$ with lowest degree term
$$r^{-(k+1)} \cdot \pi_*(\xi^{k+1}) = r^{-(k+1)} \eta_{0,k+1}\, .$$ 

Since the second term in Lemma~\ref{Lem:GRR} is supported on $\Delta$ before
$\pi$ push-forward, we can write the
result after $\pi$ push-forward as a sum of contributions of
graphs $\Gamma$ with a single edge $e=(h,h')$. Since 
the factor $e^{\xi/r}$ is then the $\pi$ pull-back  of $\xi_e$, the crucial
calculation is 
\begin{equation}
\label{ccch}
\pi_*\left[ e^{-D/r} \left( 1 + j_*\, \frac{1}{2} \sum_{m \geq 1} \frac{B_{2m}}{(2m)!} \frac{\nu_1^{2m-1} + \nu_2^{2m-1}}{\nu_1 + \nu_2}
\right)
\right]\, .
\end{equation}

Fortunately, the $\pi$ push-forward \eqref{ccch}
 was computed by Chiodo in Step 3 of \cite[Section 3]{Chiodo2}. 
The codimension $k$ part of \eqref{ccch}
is equal to a sum over $r$-weighted prestable decorated graphs $\Gamma$
with a single edge $e=(h,h')$ with coefficient
\begin{equation}\label{ll2399}
r \, \frac{B_{k+1}(\tw(h)/r)}{(k+1)!}
\sum_{i+j = k-1} (-1)^i \psi_h^i \psi_{h'}^j\, .
\end{equation}
Here, $B_k(x)$ is the Bernoulli polynomial defined by
$$
\frac{t \, e^{xt}}{e^t-1} = \sum_{k \geq 0}
B_k(x) \frac{t^k}{k!} \, .
$$
We see \eqref{ll2399} is a 
a polynomial in $\tw(h)$, $r$ and $r^{-1}$ with lowest degree term in~$r$
given by 
\begin{equation}\label{ll23999}
r^{-k} \, \frac{\tw(h)^{k+1}}{(k+1)!}
\sum_{i+j = k-1} (-1)^i \psi_h^i \psi_{h'}^j\, .
\end{equation}
The formula for the second term in Lemma \ref{Lem:GRR} is then 
obtained by multiplying by $e^{\xi_e/r}$.
\qed
\vspace{8pt}

\begin{proposition} \label{Prop:class}
The GRR formula for the Chern class $c_k(-R\pi_* \tL)$ is a polynomial family with valuation $-2k$. The coefficient of lowest degree in~$r$ equals the degree~$k$ part of the mixed degree class
\begin{multline*}
\hspace{-10pt}\sum_{
\substack{\Gamma\in \G_{g,n}^{r,L}} 
%\\w\in \mathsf{W}_{\Gamma,r}}
}
\frac{r^{-2k+|E(\Gamma)|}}{|\Aut(\Gamma)| }
\;
j_{\Gamma*}\Bigg[
%\prod_{i=1}^n \exp\left(\frac12 a_i^2 \psi_i + a_i \xi_i \right)
\prod_{v \in \V(\Gamma)} \exp\left(-\frac12 \eta(v) \right)
\\ \hspace{+10pt}
\prod_{e=(h,h')\in \E(\Gamma)}
\frac{1-\exp\left(-\frac{\tw(h)\tw(h')}2(\psi_h+\psi_{h'})\right)}{\psi_h + \psi_{h'}} \Bigg]\, ,
\end{multline*}
where $\eta=\eta_{0,2}$. 
\end{proposition}

\paragraph{Proof.} 
The Chern class $c_k$ can be expressed as a quasi-homogeneous polynomial in
the  Chern characters,
\begin{equation}\label{pp9911}
c_k = \frac1{k!}\ch_1^k + \dots\, .
\end{equation}

By Proposition~\ref{Prop:product}, each monomial $M$ on
the right side of \eqref{pp9911}
 is a polynomial family of classes with valuation at least
 $-(k + \deg(M))$. 
Since the monomial of largest degree is $$M = \frac{1}{k!}\ch_1^k\, ,$$
the Chern class $c_k$ is a polynomial family with valuation at least $-2k$. 
Since only the lowest degree coefficient in $r$ of the polynomial class $\ch_1$ contributes to the lowest degree in $r$ of the polynomial family~$c_k$, 
the valuation of the latter is exactly $-2k$.

More precisely, the lowest degree coefficient in 
$$1 + c_1(-R\pi_* \tL) + c_2(-R\pi_* \tL) + \dots$$ is obtained by exponentiating the formula for $\ch_1(-R\pi_*\tL)$ given by
Proposition~\ref{Prop:character} after changing the sign. A parallel
exponentiation is taken in \cite[Section 1]{JPPZ}.
\qed
\vspace{8pt}

We will now push-forward the GRR formula of Proposition \ref{Prop:class}
along the morphism
$$\epsilon:\ooMM_{g,n}^{r,L} \to \ooMM_{g,n}^Z\, $$
 from the commutative diagram~\eqref{Eq:CommDiag}.

\begin{corollary} \label{Cor:pushforwardck}
The push-forward $\epsilon_*c_k(-R\pi_* \tL)$ is a Laurent polynomial in~$r$ with valuation $2g-2k-1$ for $r$ sufficiently large. 
The coefficient of $r^{2g-2k-1}$ is obtained by substituting
 $r=0$ into the polynomial
\begin{multline*}
  \sum_{\Gamma\in \G_{g,n}^{Z}}
  \sum_{r{\text{\em -twist}\, } \tw }
 \frac{r^{-h^1(\Gamma)}}{|\Aut(\Gamma)|}
\;
j_{\Gamma*}\Bigg[
%\prod_{i=1}^n \exp\left(\frac12 a_i^2 \psi_i + a_i \xi_i \right)
\prod_{v \in \V(\Gamma)} \exp\left(-\frac12 \eta(v) \right)
\\ \hspace{+10pt}
\prod_{e=(h,h')\in \E(\Gamma)}
\frac{1-\exp\left(-\frac{w(h)w(h')}2(\psi_h+\psi_{h'})\right)}{\psi_h + \psi_{h'}} \Bigg]
\end{multline*} 
and extracting the part of degree~$k$.
\end{corollary}

\paragraph{Proof.} 
The claim follows directly from Proposition
\ref{Prop:pushforward}, Proposition \ref{Prop:class}, and the
calculation:
$$-2k+|\E(\Gamma)| +2g-2h^1(\Gamma)-|\V(\Gamma)| - (2g-2k-1) = -h^1(\Gamma)\ . $$
While graphs $\Gamma\in \G_{g,n}^{r,L}$
in Proposition \ref{Prop:class} correspond to classes
in the Chow cohomology of $\ooMM_{g,n}^{r,L}$,  graphs
$$\Gamma\in \G_{g,n}^{Z}$$
here correspond to classes
in the Chow cohomology of $\ooMM_{g,n}^Z$.
\qed

Finally, as proven by Chiodo \cite[Proposition 4.3.3]{Chiodo} and
\cite[Lemma 2.2.5]{Chiodo2}, all of the following three
push-forwards yield the same complex on $\ooMM_{g,n}^{r,L}$ :
\begin{itemize}
\item $R\pi_* \cL$ via $\pi: \ooCC_{g,n}^{r,L} \to \ooMM_{g,n}^{r,L}$\, ,
\item $R\pi_*{\coarseL}$ via $\pi: \mathsf{C}_{g,n}^{r,L} \to \ooMM_{g,n}^{r,L}$\, ,
\item $R\pi_* \tL$ via $\pi: \tC_{g,n}^{r,L} \to \ooMM_{g,n}^{r,L}$\, .
\end{itemize}
Hence, Corollary \ref{Cor:pushforwardck} also holds for $R\pi_* \cL$.

\begin{corollary} \label{Cor:pushforwardck2}
The push-forward $\epsilon_*c_k(-R\pi_* \cL)$ is a Laurent polynomial in~$r$ with valuation $2g-2k-1$ for $r$ sufficiently large. 
The coefficient of $r^{2g-2k-1}$ is obtained by substituting
 $r=0$ into the polynomial
\begin{multline*}
  \sum_{\Gamma\in \G_{g,n}^{Z}}
  \sum_{r{\text{\em -twist}\, } \tw }
 \frac{r^{-h^1(\Gamma)}}{|\Aut(\Gamma)|}
\;
j_{\Gamma*}\Bigg[
%\prod_{i=1}^n \exp\left(\frac12 a_i^2 \psi_i + a_i \xi_i \right)
\prod_{v \in \V(\Gamma)} \exp\left(-\frac12 \eta(v) \right)
\\ \hspace{+10pt}
\prod_{e=(h,h')\in \E(\Gamma)}
\frac{1-\exp\left(-\frac{w(h)w(h')}2(\psi_h+\psi_{h'})\right)}{\psi_h + \psi_{h'}} \Bigg]
\end{multline*} 
and extracting the part of degree~$k$.
\end{corollary}

\section{Localization analysis}

\subsection{Overview}
Let $X$ be a nonsingular projective variety over $\C$. 
Let 
$ S \rightarrow X $
be a line bundle.
Let 
$$A=(a_1,\ldots,a_n)$$
be 
a vector of double ramification data as defined in Section~\ref{Ssec:RubberToX},
$$\sum_{i=1}^n a_i= \int_\beta c_1(S)\, .$$
The double ramification cycle
$$\DR_{g,A,\beta}(X) \in A_*(\oM_{g,n,\beta}(X))$$
is defined via the moduli space of stable maps to rubber
$\oM_{g,A,\beta}^\sim(X,S)$.
We prove here the claim of Theorem~\ref{Thm:main},
 $$\DR_{g,A,\beta}(X) =  \mathsf{P}_{g,A,\beta}^g \in A_*(\oM_{g,n,\beta}(X))\, .$$
Our path follows \cite[Section 2]{JPPZ}.

  Theorem \ref{Thm:main}
for an arbitrary vector $A$ can be deduced from the case where
every $a_i$ is nonzero
by forgetting the markings with $a_i=0$. Our proof of Theorem \ref{Thm:main}
has no mathematical difficulties when $A$ has zeros, but the discussion
then requires separating markings into three types instead of
just positive and negative. For simplicity, we assume every $a_i$
is nonzero.

\subsection{Target geometry}
Denote by $\PP(X,S)$ the $\CP^1$-bundle
$\PP(\cO_X \oplus S) \rightarrow X$.
Following the notation of Section \ref{Ssec:RubberToX}, let  
$$ D_0\, ,\, D_\infty \, \subset \, \PP(X,S)$$
be the divisors defined by the projectivizations of
the loci $$\cO_X \oplus \{0\}\ \ \text{and} \ \ \{0 \} \oplus S$$
respectively.
%We will call $D_0$ the $0$-divisor and 
% $D_\infty$ the $\infty$-divisor.
After applying Cadman's $r$th root construction \cite{Cadman} to
$D_0$, we obtain a bundle
$$\PP(X,S)[r] \to X$$
with fiber given by the orbifold projective line  $\CP^1[r]$ with single
orbifold point with stabilizer $\Z/r\Z$.

Denote by $\oM_{g, A, \beta}(\PP(X,S)[r]/D_\infty)$
the moduli space of stable maps to the orbifold $\PP(X,S)[r]$ relative
to $D_\infty$. 
%Denote by $n_0$ the number of elements of $A$ equal to $0$. Then
The moduli space parametrizes connected, nodal, twisted curves~$(C,x_1,\ldots, x_n)$ of genus~$g$
with $n$ markings{\footnote{As always, the markings are
    distinct and away from the nodes.}}
together with a map 
$$f: C \to P(X,S)[r]\, ,$$ 
where $P(X,S)[r]$ is an expansion{\footnote{The expansion has a canonical $D_0\subset
    P(X,S)[r]$ from the bulk and a canonical \newline $D_\infty \subset P(X,S)[r]$ from the last component
    of the expansion.}} 
of $\PP(X,S)[r]$ along $\D_\infty$. 
The following conditions are required to hold over
$$D_0,D_\infty\subset P(X,S)[r]\, :$$
\begin{enumerate}
\item[(i)] The stack structure of the domain curve $C$ occurs
only at the nodes over $D_0$
  and at the markings corresponding to positive elements of~$A$ (which must be mapped to $D_0$). The monodromies associated to the latter
  markings are specified by the parts $a_i$ of~$A$ at these markings. More precisely,
  the monodromies are $a_i \bmod r$.

\item[(ii)] The only points mapped by $f$ to $D_\infty$
are the markings $x_i$ with negative $a_i$. The multiplicity of $D_\infty$ at
  $x_i$ is $-a_i$. 
The map $f$ satisfies the
ramification matching condition over the internal nodes of the expansion~$P$.
\item[(iii)] Stability requires the full data $(C,x_1,\ldots,x_n, f, P(X,S)[r])$
  to have only finitely many automorphisms.
\end{enumerate}

The moduli space $\oM_{g, A,\beta}(\PP(X,S)[r]/D_\infty)$
has a perfect obstruction theory and a virtual class of dimension
\begin{multline}\label{zzzz1}
\dim_{\C}\, [\oM_{g, A, \beta}(\PP(X,S)[r]/D_\infty)]^{\vir}= \\
\dim_\C\, [\oM_{g,n,\beta}(X,S)]^\vir - (g-1) + \sum_{i | a_i >0} \lfloor \frac{a_i}r \rfloor\, ,
\end{multline}
see \cite[Section 1.1]{JPT}.

We will be most interested in the  case where $r> \sum_{i=1}^n |a_i|$. The last term in
\eqref{zzzz1} then vanishes. We refer the reader to \cite{AGV, JPT,MauPan} for a more detailed definition of the moduli space of stable relative maps.

\subsection{The $\C^*$-fixed loci}
\subsubsection{The $\C^*$-action}
The standard $\C^*$-action over $X$
on the projective bundle $\PP(\cO \oplus S) \to X$, defined by
$$\xi\cdot[x,z,s] =[x,z,\xi s]\, ,$$
lifts canonically to $\C^*$-actions on
$\PP(X,S)[r]$ and $\oM_{g,A,\beta}(\PP(X,S)[r]/D_\infty)$.
%
%
%We want to consider virtual localization for $\oM_{g, n,
%  \mu}(\PP[r], \nu)$ for the $\C^*$-action on $\PP[r]$ with weights
%$(0, 1)$. Let $t$ be the first Chern class of $\C$ with weight one
%$\C^*$-action. 
%We describe here the $\C^*$-fixed loci of 
 %$\oM_{g,A,\beta}(\PP(X,S)[r],D_\infty)$.

\subsubsection{Graphs}\label{ffpp33}
The $\C^*$-fixed loci of $\oM_{g,A,\beta}(\PP(X,S)[r]/D_\infty)$
are in bijective correspondence with decorated graphs~
$$\Phi=(\V\, ,\ \E\, ,\ \L\, , \ \mathrm{g}:\V \rarr \Z_{\geq 0}\, ,\ 
\beta: \V \to H_2(X,\Z)\, ,\ \ell:\V\rarr \{0,\infty\}\, ,
\ d:\E \rarr \Z_{>0}) \, $$
which satisfy the following six properties:
 
%\begin{enumerate}
%\item[(i)] $\V$ is a vertex set with a genus function $\g:V\to \Z_{\geq 0}$,
%\item[(ii)] $\H$ is a half-edge set equipped with a 
%vertex assignment $v:\H \to \V$ and an involution $\iota$,
%\item[(iii)] $\E$, the edge set, is defined by the
%2-cycles of $\iota$ in $\H$ (self-edges at vertices
%are permitted),
%\item[(iv)] $\L$, the set of legs, is defined by the fixed points of $\iota$ an%d is
%placed in bijective correspondence with a set of $n$ markings,
%\item[(v)] the pair $(\V,\E)$ defines a {\em connected} graph
%satisfying the genus condition 
%$$\nice \sum_{v \in \V} \g(v) + h^1(\Gamma) = g\, ,$$
%\item[(vi)] for each vertex $v$, the stability condition holds: if $\beta(v)=0$%,
%then
%$$2\g(v)-2+ \n(v) >0\, ,$$
%where $\n(v)$ is the valence of $\Gamma$ at $v$ including 
%both half-edges and legs,
%\item[(vii)] the
%degree condition holds:
%$$\nice \sum_{v \in \V} \beta(v) = \beta\, .$$
%\end{enumerate}

\begin{enumerate}
\item[(i)] $\V$ is a vertex set with a genus function $\g$, a degree function
$\beta$, and a label $\ell$.
For $v\in \V$, the degree $\beta(v)$ must be an effective{\footnote{Effective
here includes the class $0\in H^2(X,\Z)$.}} curve
class. We also require the genus and degree  conditions to hold: 
$$ g=\sum_{v\in \V} \g(v) + h^1(\Phi) 
\ \ \ \ \text{and}\ \ \ \ \beta=\sum_{v\in \V} \beta(v)\, .$$

\item[(ii)] $\L$, the set of legs, is placed in bijective correspondence
with the $n$ markings:
\begin{itemize}
\item legs marked $i$ with $a_i>0$ are incident to vertices labeled~$0$, 
\item legs marked $i$ with $a_i<0$ are incident to vertices labeled~$\infty$. 
%\item Legs marked $i$ with $a_i=0$ can be incident to vertices of either label.
\end{itemize}

\item[(iii)] $\E$ is the edge set. For 
$e\in \E$, the edge degree
 $d_e$ corresponds to the $d_e$-th power map $$\CP^1[r] \to \CP^1[r]\, .$$

\item[(iv)] $\Phi$ is a connected graph, and $\Phi$ is bipartite 
 with respect to labeling $\ell$: every edge
  is incident to a 0-labeled vertex and an $\infty$-labeled vertex.

\item[(v)] If $\ell(v) =0$, denote by $A(v)$ the list of integers formed 
by the values $a_i$ for the legs $i$ incident to $v$ {\em and} by the values
$-d_e$ for the edges $e$ incident to~$v$. 
For every such vertex $v$, we impose the condition 
$$|A(v)| = \int_{\beta(v)} c_1(S) \ \ \bmod r\, ,$$
where $|A(v)|$ is the sum of the elements of $A(v)$.
 %For $r$ large enough this becomes simply $|A(v)| = (\xi,\beta(v))$, because the numbers $(\xi, \beta')$ for an effective summand $\beta'$ of $\beta$ are bounded.

\item[(vi)] If $\ell(v)=\infty$, denote by $A(v)$ the list of integers formed by the values $a_i$ for the legs $i$ incident to $v$ {\em and} by the values
 $d_e$ for the edges $e$ incident to~$v$. For every such vertex, we impose the condition $$|A(v)| = \int_{\beta(v)} c_1(S) \ .$$

\end{enumerate}

To every 0-labeled vertex $v$ of $\Phi$, we assign the space $\oM_{\g(v), A(v),\beta(v)}^r(X,S)$. We will use the notation
$$\oM_v^r\, =\, \oM^r_{\g(v),A(v),\beta(v)}(X,S)\, .
$$
As  explained in Section~\ref{Ssec:rthroot}, the forgetful map
$$\oM_v^r\rightarrow \oM_{\g(v),\n(v), \beta(v)}(X)$$
is a finite map of degree~$r^{2g-1}$. 
The virtual fundamental class  
 $\left[ \oM_v^{r} \right]^{\vir}$ of $\oM_v^r$ is of dimension
\begin{eqnarray*}
\dim_{\C}\, \left[ \oM_v^{r} \right]^{\vir}& =&\dim_{\C}\, \left[ \oM_{\g(v), \n(v), \beta(v)}(X) \right]^{\vir}\\
& =& 
(3-\dim(X))(\g(v)-1)+ \n(v) +\int_{\beta(v)} c_1(X) \,  .
\end{eqnarray*}

A $\C^*$-fixed relative stable map with vector $A=(a_1,\ldots,a_n)$, 
$$f:(C,x_1,\ldots,x_n) \to \PP(X,S)[r]/D_\infty$$
takes two basic forms:
\begin{itemize}
\item If the target does not expand,
 then the stable map has a finite number of preimages of $D_\infty$
which correspond precisely to the markings $i$ with $a_i<0$. 
Each such preimage is then described by an unstable vertex of $\Phi$ 
decorated by~$\infty$. 
\item 
If the target expands, then $C$ contains  a  possibly 
disconnected subcurve mapping to the rubber 
 $\PP(X,S)$ -- with no orbifold structure at $D_0$ in the rubber.
The ramification in the rubber over $D_0$ is specified by the degrees $d_e$ of the edges of $\Phi$, and the ramification over $D_\infty$ is specified by the negative elements of $A$. 
The ramification in the rubber occurs over $D_0$ and is specified by
 the negative elements of $A$ over $D_\infty$  
Every $\infty$-labeled vertex $v$ describes a connected component of the rubber map. We will denote the moduli space
of stable maps to rubber by $\oM_\infty^{\sim}$.
\end{itemize}

In the second case above,
let $\g(\infty)$ be the genus of the possibly disconnected domain of the rubber map, $\n(\infty)$ the total number of legs and edges adjacent to $\infty$-labeled vertices, and $\beta(\infty)$ their total degree. The
 virtual fundamental class  $\left[\oM_\infty^{\sim} \right]^{\vir}$ has dimension
\begin{eqnarray*}
\dim_{\C}\, \left[ \oM_\infty^{\sim} \right]^{\vir}&=&
\dim_\C \left[ \oM_{\g(\infty),\n(\infty), \beta(\infty)}(X) \right]^\vir - \g(\infty) \\
&=& 
(3- \dim(X))(\g(\infty)-1)+ \n(\infty) + \int_{\beta(\infty)}
c_1(X) - g(\infty)\, \, .
\end{eqnarray*}
The image of the virtual fundamental class 
$\left[\oM_\infty^{\sim} \right]^{\vir}$
in the moduli space of (not necessarily connected) stable maps to~$X$ is denoted by $\DR_\infty$. 

\subsubsection{Unstable vertices}
A vertex $v\in\mathsf{V}(\Phi)$ is {\em unstable} if $\beta(v)=0$ and $2 \g(v) - 2 + \n(v) \leq 0$. There are four types of unstable vertices: 
\begin{enumerate}
\item[(i)] $\ell(v)= 0$, $\g(v)=0$, $v$ carries no markings and
 one incident edge,
\item[(ii)] $\ell(v)= 0$, $\g(v)=0$,  
$v$ carries no markings and two incident edges,
\item[(iii)] $\ell(v)= 0$, $\g(v)=0$, 
 $v$ carries one marking and one incident edge,
\item[(iv)] $\ell(v)=\infty$, $\g(v)=0$, 
$v$ carries one marking and one incident edge.
\end{enumerate}
The target of the stable map expands if and only if
 there is at least one $\infty$-labeled stable vertex.

A stable map in the $\C^*$-fixed locus
corresponding to  $\Phi$ is obtained by gluing together maps
associated to the vertices $v\in \mathsf{V}(\Phi)$ with
Galois covers 
associated to the edges. Denote by $\mathsf{V}_{\rm st}^0(\Phi)$ the set of 0-labeled stable vertices of~$\Phi$. Then the $\C^*$-fixed locus
corresponding to $\Phi$ is isomorphic to the product 
$$
\oM_\Phi =
 \begin{cases}
      \prod\limits_{v\in \mathsf{V}^0_{\rm st}(\Phi)}\, \oM_v^r\ \times\ \oM_\infty^{\sim}      \, , & \text{if the target expands,} 
\vspace{7pt}\\ 
    \prod\limits_{v\in \mathsf{V}^0_{\rm st}(\Phi)}\, \oM_v^{r}
\, , & \text{if the target does not expand,}
  \end{cases}
$$
quotiented by the automorphism group of $\Phi$ {\em and} the product of 
cyclic groups $\Z_{d_e}$ associated to the Galois covers of the edges.

%In case $v_\infty$ is stable, there is a multiplicity factor of $\prod_{e\in \mathsf{E}(\phi)} d_e$ from the geometry of the moduli space of stable relative maps, see \cite{GrV}.

The natural morphism corresponding to $\Phi$,
$$\iota: \oM_\Phi \rightarrow \oM_{g, A, \beta}(\PP(X,S)[r]/D_\infty)\, ,$$
is of degree
$$|\Aut(\Phi)|  \cdot \prod_{e\in \mathsf{E}(\Phi)}d_e  $$
%$$
%\begin{cases}
%|\Aut(\Phi)|
% & \text{if the target expands, }
% \vspace{7pt}\\ 
%|\Aut(\Phi)|  \cdot \prod_{e\in \mathsf{E}(\Phi)}d_e     & \text{if the target %does not expand,} 
%  \end{cases}
%$$
onto the image $\iota(\oM_{\Phi})$.

\begin{lemma}\label{ggee99}
For $r$ sufficiently large, the unstable vertices of type (i) and (ii)
can {\em not} occur.
\end{lemma}

\paragraph{Proof.} 
We define  $\beta' \in H_2(X,\Z)$ to be an
 {\em effective summand} of $\beta$ if both $\beta'$ and $\beta-\beta'$ are effective cycle classes (including 0). 
Let $b$ be the maximum of 
%$|(\beta',\xi)|$ 
$\left|\, \int_{\beta'} c_1(S) \, \right|$
over all effective summands of $\beta$.
 Further, let 
$$a_+ = \sum_{i| a_i >0} a_i$$ 
be the sum of the positive elements of the vector $A$.
% that is, the intersection index of the degree of the stable map in $\PP(X,S)[r]$ with the 0-divisor $D_0$. 
Assume  
$$r > 2(a_+ + b)\, .$$

Let $\beta_0$ (respectively, $\beta_\infty$) be the sum of degrees of all vertices of $\Phi$ with label~0 (respectively, $\infty$), so 
 $$\beta=\beta_0 + \beta_\infty\, .$$
 We then have  $$
a_+ - \int_{\beta_0}c_1(S) = \sum_{e\in \E(\Phi)} d_e\, .
$$
By our choice of $r$, we have $\sum_{e\in \E(\Phi)} d_e < r/2$.

At each 0-labeled stable vertex $v\in \mathsf{V}(\Phi)$, the condition
$$
 \sum_{i \, \vdash v} a_i - \int_{\beta(v)}c_1(S) = \sum_{e \, \vdash v} d_e\ \bmod r
$$
holds by the conditions on the graph~$\Phi$. By our choice of $r$, both the absolute value of $\sum_{i \, \vdash v} a_i - \int_{\beta(v)}c_1(S)$ and 
$\sum_{e \, \vdash v} d_e$ are less than $r/2$.
Therefore, the equality mod~$r$ is actually an exact equality:
\begin{equation}
\label{llqq88}
\sum_{i \, \vdash v} a_i - \int_{\beta(v)}c_1(S) = \sum_{e \, \vdash v} d_e\, . 
\end{equation}

For an unstable vertex of type~(i) or~(ii), we have both $\sum_{i \, \vdash v} a_i = 0$ and $\beta(v)=0$. The sum of the degrees of edges adjacent to such a vertex then vanishes by \eqref{llqq88}. However, the degree of every edge is a positive integer and the graph~$\Phi$ is connected. 
The resulting contradiction implies that there are 
no unstable vertices of types~(i) and~(ii).
 \qed

\subsection{Localization formula}
\label{locformula}

We write the $\C^*$-equivariant Chow ring of a point
as $$A^*_{\C^*}(\bullet) = \Q[t]\, ,$$
where $t$ is the first Chern class of the standard representation.

For the localization formula, we will
require the inverse of the $\C^*$-equivariant Euler
class of the virtual normal bundle in $\oM_{g,A,\beta}(\PP(X,S)[r]/D_\infty)$
to the $\C^*$-fixed
locus corresponding to~$\Phi$. Let
$$f: (C,x_1,\ldots,x_n) \to P(X,S)[r]\, , \ \ \ \ \ [f]\in \oM_\Phi\, ,$$
where $P(X,S)[r]$ is a possible expansion of $\PP(X,S)[r]$ along $D_\infty$.
Denote by $$T \to \PP(X,S)[r]$$
 the tangent line bundle to the fiber of $\PP(X,S)[r] \to X$. 
%It is easy to see that $T(-D_\infty)|_{D_\infty}$ is trivial.
For simplicity, we 
will also denote by $T$ the pull-back of $T$ from $\PP(X,S)[r]$ to the
 expansion $P(X,S)[r]$.
The formula for the inverse Euler class can then
be written as: 
\begin{equation}
  \label{eq:contribs}
\frac{1}{e(\text{Norm}^{\vir})}=  \frac{e(H^1(C, f^*T(-D_\infty)))}{e(H^0(C, f^*T(-D_\infty)))} \frac 1{\prod_i e(N_i)} \frac 1{e(N_\infty)}\, .
\end{equation}

Several aspects of 
Formula \eqref{eq:contribs} require
explanation. 
To start, we assume  that~$r$ is sufficiently large (using Lemma \ref{ggee99})  
to exclude the presence of
 unstable vertices of types~(i) and~(ii) in~$\Phi$. 
To compute the leading factor of~\eqref{eq:contribs},
\begin{equation}\label{xhh}
\frac{e(H^1(C, f^*T(-D_\infty)))}
{e(H^0(C, f^*T(-D_\infty)))} \, ,
\end{equation}
we use the normalization exact sequence for the domain $C$ tensored with the line bundle $f^*T(-D_\infty)$. 
The associated long exact sequence in cohomology decomposes the 
leading factor into a product of vertex, edge, and node contributions:

\begin{enumerate}

\item[$\bullet$]
Let $v\in \mathsf{V}(\Phi)$ be a stable vertex over $D_0\subset \PP(X,S)[r]$ 
corresponding to 
a moduli space $$\oM_v^r = \oM^r_{\g(v), A(v), \beta(v)}(X,S)\, .$$
 The orbifold universal curve{\footnote{
The moduli space $\oM^r_{g,A,\beta}(X)$
may be considered with the universal curve
$\mathcal{C}^r_{g,A,\beta}$ of Section \ref{www111}
or with the orbifold universal curve
$\mathcal{C}^{r,\text{orb}}_{g,A,\beta}$.
While $\mathcal{C}^r_{g,A,\beta}$ 
has orbifold structure only at the nodes of the fibers,
$\mathcal{C}^{r,\text{orb}}_{g,A,\beta}$ 
has orbifold structure {both}  at the markings $x_i$ {\em and}
at the nodes. A full discussion of the differences
will be given in Section \ref{asdasd}.}} 
$$\pi:\mathcal{C}^{r,\text{orb}}_{\g(v),A(v),\beta(v)} \rightarrow \oM^r_{\g(v),A(v),\beta(v)}$$
 carries an orbifold line bundle $\cL^{\text{orb}}$ (the $r$th root of the pull-back of $S$)  which is the
 pull-back of $T$ to the universal curve. Therefore,  the
contribution 
\begin{equation*}
  \frac{e(H^1(C_v, f^*T(-D_\infty))}
  {e(H^0(C_v, f^*T(-D_\infty))}
\end{equation*}
yields the class
\begin{equation*}  
e\left((-R{\pi}_*\cL^{\text{orb}}) \otimes \cO^{(1/r)}\right) \
=
c_{\mathrm{rk}}\left((-R{\pi}_*\cL^{\text{orb}}) \otimes \cO^{(1/r)}\right) \
\end{equation*}
in $A^*(\oM^r_v) \otimes \Q\left[t,\frac 1t\right]$,
where $\cO^{(1/r)}$ is a trivial
line bundle with a $\C^*$-action of weight $\frac 1r$ and
$$\mathrm{rk} = \g(v)-1 + |\mathsf{E}(v)|$$
is the virtual rank of $-R\pi_*\cL^{\text{orb}}$.

Unstable $0$-labeled vertices of type~(iii) contribute factors of~1. 
Since the restriction of $T(-D_\infty)$ to $D_\infty$ is trivial,
the $\infty$-labeled vertices also contribute factors of~1.

\item[$\bullet$] The edge contribution is trivial since the degree $\frac {d_i}r$ of
$f^*T(-D_\infty)$ is less than 1, see \cite[Section 2.2]{JPT}.

\item[$\bullet$] The contribution of a node $N$ over $D_0$ is trivial. Indeed, the
space of sections
 $H^0(N,f^*T(-D_\infty))$ vanishes  because $N$ must be
stacky, and $H^1(N,f^*T(-D_\infty))$ is trivial for dimension reasons.
Nodes over $D_\infty$ contribute~1.
\end{enumerate}

\noindent Consider next the last two factors of \eqref{eq:contribs},
$$\frac 1{\prod_i e(N_i)} \frac 1{e(N_\infty)}\, .$$

\begin{enumerate}
\item[$\bullet$]
The  product ${\prod_i e(N_i)^{-1}}$
is over the nodes that correspond to half-edges of the graph~$\Phi$ adjacent to a $0$-labeled vertex.
If $N$ is a node corresponding to an edge $e\in \mathsf{E}(\Phi)$ 
and the associated vertex~$v$ is stable, then 
\begin{equation} \label{Eq:eN}
 e(N) = \frac {t+\text{ev}_e^*(c_1(S))}{r \, d_e} - \frac{\psi_e}r\, .
\end{equation}
The factor corresponds
 to the smoothing of the node $N$ of the domain curve: $e(N)$ is 
the first Chern class of the normal line bundle of 
the divisor of nodal domain curves.
The first Chern classes of the tangent lines to the branches at the node 
are divided by $r$ because of the orbifold twist, see Section
\ref{Ssec:twistedcurves}.

In the case of an unstable vertex of type (iii), the associated
edge does {\em not}
produce a node of the domain. The type (iii) edge incidences do
{\em not} appear in ${\prod_i e(N_i)^{-1}}$.
%$N$ is {\em not} a node over $0$ and the associated vertex is unstable of type %(iii), then
%\begin{equation*}
%  e(N) = \frac 1{d_e}.
%\end{equation*}

\item[$\bullet$]
$N_\infty$ corresponds to the expansion of the target $\PP(X,S)[r]$ over $D_\infty$.
The factor $e(N_\infty)$
is 1 if the target $(\PP(X,S)[r]/D_\infty)$ does not expand and 
$$e(N_\infty)= -\frac{t +\Psi_\infty}{\prod_{e\in \mathsf{E}(\Phi)}d_e} $$
 if the target expands. 

Here, $\Psi_\infty$ is the first Chern class of a line bundle defined as follows. Consider a point of the moduli space $\oM_{g,A,\beta}(\PP(X,S)[r]/D_\infty)$ where the target expands. For the target over the point, 
the divisor along which the target expands carries tangent line bundles to the two components of the target. The tensor product of these two line bundles is a trivial line bundle. Thus the tensor product is the pull-back of a line bundle $N$ over the divisor of $\oM_{g,A,\beta}(\PP(X,S)[r]/D_\infty)$ where the target expands, and
 $\Psi_\infty$ is the first Chern class of $N$.
See~\cite{MauPan} for more details.
\end{enumerate}

The virtual class
of $\oM_{g,A,\beta}(\PP(X,S)[r]/D_\infty)$ can be written in terms of
the $\C^*$-fixed point loci by the virtual localization formula~\cite{GrP}:
\begin{multline}\label{llff}
\left[\oM_{g, A, \beta}(\PP(X,S)[r]/D_\infty) \right]^{\vir} =\\
\sum_{\Phi} \frac{1}{|\Aut(\Phi)|} \, \frac{1}
{\prod_{e\in\mathsf{E}(\Phi)}d_e}\cdot 
\iota_*\left(\frac{[\oM_\Phi]^{\vir}}
{e(\text{Norm}^{\vir})}\right)\,
\end{multline}
in $A^*\left(\oM_{g, A, \beta}(\PP(X,S)[r]/D_\infty)\right)\otimes 
\Q[t,\frac{1}{t}]$.
Our analysis of the inverse Euler class of the virtual normal bundle
yields the following contributions to $\frac{[\oM_\Phi]^{\vir}}
{e(\text{Norm}^{\vir})}$ associated to the graph $\Phi$ :
\begin{itemize}
\item a factor
\begin{equation*}
\prod_{e \in \mathsf{E}(v)} \frac {r}{\frac{t+\text{ev}_e^*(c_1(S))}{d_e}-\psi_e}\, \cdot\, 
\sum_{d \geq 0} c_d(-R\pi_*\cL^{\text{orb}}) \left(\frac{t}r\right)^{\g(v)-1+|\E(v)|-d}
\end{equation*}
for each stable vertex $v\in \mathsf{V}(\Phi)$ over 0, where $\E(v)$ is the
set of edges incident to $v$,
\item a factor 
\begin{equation*}
   -\frac {{\prod_{e\in\mathsf{E}(\Phi)}d_e}       }{t+ \psi_\infty}\cdot {\mathsf{DR}_\infty} 
\end{equation*}
if the target expands, where $\mathsf{DR}_\infty$ is the virtual class
of the moduli space of map to the rubber over $D_\infty$.
%(which is equivalent to $v_\infty$ being stable).
\end{itemize}

\subsection{The formula for the DR-cycle}

\subsubsection{Three operations}
We will now perform three operations on the localization formula \eqref{llff}
for the virtual class $[\oM_{g, A, \beta}(\PP(X,S)[r]/ D_\infty)]^{\vir}$:

\begin{enumerate}
\item[(i)] the $\C^*$-equivariant push-forward via
\begin{equation}\label{xx88}
\epsilon: \oM_{g, A, \beta}(\PP(X,S)[r]/ D_\infty) \rightarrow \oM_{g,n,\beta}(X)
\end{equation}
to the moduli space $\oM_{g,n,\beta}(X)$ of stable maps to~$X$ with trivial $\C^*$-action, 
\item[(ii)] extraction of the coefficient of $t^{-1}$ after push-forward by $\epsilon_*$, 
\item[(iii)] extraction of the coefficient of $r^0$.
\end{enumerate}

After push-forward by $\epsilon_*$, the coefficient of $t^{-1}$
is equal to 0 because 
$$\epsilon_*[\oM_{g, A, \beta}(\PP(X,S)[r]/ D_\infty)]^{\vir} \in A_*(\oM_{g,n,\beta}(X))\otimes
\Q[t]\, .$$
Using Proposition~\ref{Prop:class}, all terms of the $t^{-1}$ coefficient will be seen to be %Laurent 
polynomials in~$r$,
so operation~(iii) will be well-defined. 
After operations~(i-iii), only two nonzero terms will remain.
The cancellation of the two remaining terms will prove Theorem~\ref{Thm:main}.

To perform (i-iii),
we multiply the $\epsilon$-push-forward of the
localization formula \eqref{llff} by $t$
 and extract the coefficient of $t^0 r^0$.
To simplify the computations, we introduce the new variable 
$$s=tr\, .$$ 
%so that $t= s/r$. 
Then, instead extracting the coefficient of $t^0r^0$,
we extract the coefficient of $s^0 r^0$. 

\subsubsection{Push-forward to $\oM_{g,n, \beta}(X)$} \label{pppfff}

For each vertex $v\in \mathsf{V}(\Phi)$, 
following the notation of Section \ref{ffpp33},
we have
$$\oM_v^r =\oM^r_{\g(v),A(v),\beta(v)}\, .$$
%DZ: next sentence rewritten.
As in diagram \eqref{Eq:CommDiag}, we denote by 
\begin{equation} \label{xx99}
\epsilon: \oM^r_v \rightarrow \oM_{\g(v),\n(v),\beta(v)}(X)
\end{equation}
the morphism obtained by forgetting the $r$th root line bundle.
The maps $\epsilon$ in \eqref{xx88} and \eqref{xx99}
are compatible. Denote by 
\begin{equation}\label{ffrrtt}
\widehat{c}_d \, = \, r^{2d-2\g(v)+1} \epsilon_* c_d(-R\pi_*\cL^{\text{orb}})
\, \in\, A^d(\oM_{\g(v),\n(v),\beta(v)}(X))\, .
\end{equation}
By Corollary~\ref{Cor:pushforwardck}, $\widehat{c}_d$  is a polynomial in $r$ for $r$ sufficiently large.

We now write
the inverse Euler class of the virtual normal bundle
for the $\C^*$-fixed point locus associated to the
graph  $\Phi$
{\em after} push-forward along 
$$\epsilon: \oM_{g,n,\beta}(\PP(X,S)[r]/D_\infty) \rightarrow \oM_{g,n,\beta}(X)$$
in terms of $s=rt$.
The analysis of Section~\ref{locformula} yields
the following contributions to $\epsilon_*\iota_*\frac{[\oM_\Phi]^{\vir}}
{e(\text{Norm}^{\vir})}$ :
\begin{itemize}
\item a factor
\begin{equation*}
  \frac{r}{s}\cdot \prod_{e \in \mathsf{E}(v)} \frac {d_e}{1+\frac{r}{s} 
\text{ev}_e^*(c_1(S))-\frac{rd_e}{s}\psi_e}\, \cdot\,
  \sum_{d \geq 0} \widehat{c}_d \, s^{\g(v)-d} \, \cdot \, 
  \left[ \oM_{\g(v),\n(v), \beta(v)}(X) \right]^\vir
  \end{equation*}
  \begin{equation*}
  \in A_*(\oM_{\g(v),\n(v), \beta(v)}(X))\otimes \Q\left[s,\frac 1s\right]
\end{equation*}
for each stable vertex $v\in \mathsf{V}(\Phi)$ over 0,
\item a factor 
\begin{equation*}
   -\frac{r}{s}\cdot \frac {{\prod_{e\in\mathsf{E}(\Phi)}d_e}       }{1+\frac{r}{s} \Psi_\infty}\cdot \mathsf{DR}_\infty 
\end{equation*}
if the target degenerates.
% (which is equivalent to $v_\infty$ being stable).
\end{itemize}

\noindent For the first factor, we have used the compatibility of
the virtual classes 
$$\epsilon_* \left[ \oM^r_v\right]^\vir = r^{2g-1}
\left[\oM_{\g(v),\n(v),\beta(v)}(X)\right]^\vir\, $$
proven in \cite[Theorem 6.8]{AJT}.

%\begin{itemize}
%\item A global factor of $\frac{s}r \frac 1{|\Aut\Gamma|}$.
%\item A factor
%\begin{equation*}
%[\oM_v] \cdot \frac rs 
%\prod_{e \in E_v} \frac {d_e}{1 - \frac{r d_e}s \psi_e} \sum_{d \geq 0} {\hat% c_d} s^{g_v-d}
%\end{equation*}
%for each stable vertex over 0.
%\item A factor 
%\begin{align*}
%   - \frac rs \frac 1{1+ \frac rs \psi_\infty} \cdot \DR_\infty^o \qquad & \m%box{if the target degenerates or}\\
%  \frac 1{\prod d_e} \qquad &\mbox{if it does not.}
%\end{align*}
%\end{itemize}
%Here $\DR_\infty^o$ stands for the possibly disconnected DR-cycle and include%s factors $\frac 1{d_e}$ for type 3 unstable vertices at $\infty$. 

\subsubsection{Extracting coefficients}

From \eqref{llff} and the contribution calculus for $\Phi$ 
presented in Section \ref{pppfff}, we have a complete formula
for the $\C^*$-equivariant push-forward of
$t$ times the virtual class:
\begin{multline}\label{llfff}
\epsilon_*\Big(t[\oM_{g, A, \beta}(\PP(X,S)[r], D_\infty)]^{\vir}\Big) = \\
\frac{s}{r}\cdot\sum_{\Phi} \frac{1}{|\Aut(\Phi)|} \, \frac{1}
{\prod_{e\in\mathsf{E}(\Phi)}d_e}\cdot 
\epsilon_*\left(\frac{[\oM_\Phi]^{\vir}}
{e(\text{Norm}^{\vir})}\right)\, .
\end{multline}

\paragraph{Extracting the coefficient of $r^0$.} 
By Corollary~\ref{Cor:pushforwardck}, 
the classes $\widehat{c}_d$ are polynomial in $r$ for $r$ sufficiently large.
%$$r > 2|\mu|\, .$$
We have an $r$ in the denominator in the prefactor
on the right side of \eqref{llfff} which
 comes from the multiplication by $t$ on the left side.
However, in all other factors, we only have positive powers of $r$,
 with at least one $r$ per $0$-labeled vertex of the graph
and one more~$r$ if the target degenerates. The {\em only} graphs $\Phi$ which
 contribute to the coefficient of $r^0$ are those with 
{\em exactly one $r$ in the numerator}. There are only two graphs which have exactly one $r$ factor in the numerator:
\begin{enumerate}
\item[$\bullet$]
the graph $\Phi'$ with a $0$-labeled stable vertex of full genus~$g$ and 
an $\infty$-labeled unstable vertex of type~(iv) for each negative element of $A$,
\item[$\bullet$]the graph $\Phi''$ with a stable $\infty$-labeled vertex of full genus~$g$ and a $0$-labeled unstable type~(iii) vertex for each positive element of~$A$. 
\end{enumerate}

No terms involving $\text{ev}^*(c_1(S))$, $\psi$ or $\Psi_\infty$ classes contribute
to the $r^0$ coefficient of either $\Phi'$ or $\Phi''$
since every $\psi$ class in the 
localization formula comes with an extra factor of $r$. 
%DZ: added sentence.
The homology class associated with $\Phi''$ is, by definition, the double ramification
cycle $\DR_{g,A,\beta}(X)$.
We can now write the $r^0$ coefficient of the right side of \eqref{llfff} as
\begin{equation}\label{llffff}
|\Aut| \cdot\text{Coeff}_{r^0}\left[\epsilon_*\Big(t[\oM_{g, A, \beta}(\PP(X,S)[r], D_\infty)]^{\vir}\Big)\right] \ =\ 
\end{equation}
\begin{equation*}
 \text{Coeff}_{r^0}\left\{\sum_{d \geq 0} \widehat{c}_d \, s^{g-d}
\, \cdot \, \left[ \oM_{g,n,\beta}(X) \right]^\vir
\right\} \ -\ 
   \DR_{g,A,\beta}(X)\, 
 \end{equation*}
in $A_*(\oM_{g,n,\beta}(X))\otimes \Q[s,\frac{1}{s}]$.
Here, $|\Aut|= |\Aut(\Phi')|=|\Aut(\Phi'')|=1$.

\paragraph{Extracting the coefficient of $s^0$.} 
The remaining powers of $s$ in \eqref{llffff} appear only
 in the classes $\widehat{c}_d$ (in the contribution of the graph $\Phi'$).
In order to obtain $s^0$, we must take $d=g$,
\begin{equation}\label{llfffff}
|\Aut|\cdot \text{Coeff}_{s^0r^0}\left[\epsilon_*\Big(t[\oM_{g, A, \beta}(\PP(X,S)[r], D_\infty)]^{\vir}\Big)\right] \ =\
\end{equation}
\begin{equation*}
\text{Coeff}_{r^0}\left\{\widehat{c}_g \, \cdot \, \left[ \oM_{g,n,\beta}(X) \right]^\vir \right\}  \ -\ 
   \DR_{g,A,\beta}(X)\, 
 \end{equation*}
in $A_*(\oM_{g,n,\beta}(X)$.

\subsubsection{Proof of Theorem \ref{Thm:main}}
\label{asdasd}
Since $\text{Coeff}_{s^0r^0}\left[\epsilon_*\Big(t[\oM_{g, A, \beta}(\PP(X,S)[r], D_\infty)]^{\vir}\Big)\right]$ vanishes, we can rewrite equality
\eqref{llfffff} as
\begin{multline}\label{jjqq}
   \DR_{g,A,\beta}(X) = \\ \text{Coeff}_{r^0}
\left[
r \epsilon_* c_g(-R\pi_*\cL^{\text{orb}})\cdot \left[ \oM_{g,n,\beta}(X) \right]^\vir
\right]
\, \in\, A_{*}(\oM_{g,n,\beta}(X))\, .
 \end{multline}
Here, we have used the definition \eqref{ffrrtt} of $\widehat{c}_g$,
 $$\widehat{c}_g = r \epsilon_* c_g(-R\pi_*\cL^{\text{orb}})\, .$$

What is the relationship between the line bundle
$$\cL \ \ \text{on}\ \ \pi:\mathcal{C}^r_{g,A,\beta}\rightarrow   \oM^r_{g,A,\beta}(X)$$
considered in Section \ref{www111} and the line bundle
$$\cL^{\text{orb}} \ \ \text{on} \ \ \pi:\mathcal{C}^{r,\text{orb}}_{g,A,\beta}
\rightarrow  \oM^r_{g,A,\beta}(X)$$
which appears in \eqref{jjqq} here? The definitions are slightly
different:
\begin{itemize}
\item $\mathcal{C}^r_{g,A,\beta}$ 
has orbifold structure only at the nodes of the fibers,
\item $\cL^{\otimes r} = f^*S(-\sum_{i=1}^n a_i x_i)$,
\item $\mathcal{C}^{r,\text{orb}}_{g,A,\beta}$ 
has orbifold structure {both} at the markings $x_i$ {\em and}
at the nodes,
\item $\cL^{\text{orb}} = f^*T$ .
\end{itemize}
The universal curve $\mathcal{C}^r_{g,A,\beta}$ is the coarsification
along the markings $x_i$ of $\mathcal{C}^{r,\text{orb}}_{g,A,\beta}$.
By considering the sheaf of
invariant sections of $\cL^{\text{orb}}$ on  $\mathcal{C}^r_{g,A,\beta}$,
we obtain an $r$th root of
$$f^*S\Big(-\sum_{i|a_i>0} a_i x_i - \sum_{e\in \E} (r-d_e)x_e\Big) = 
f^*S\Big(-\sum_{i=1}^n a_ix_i\Big) \otimes
{\mathcal{O}_C}\Big(-\sum_{i|a_i<0} x_i\Big)^{\otimes r}\, .$$
So the $r$th roots corresponding to $\cL$ and 
the coarsification of $\cL^{\text{orb}}$ are related simply by
the factor ${\mathcal{O}_C}\Big(-\sum_{i|a_i<0} x_i\Big)$ which
yields a shift
$$ \xi \ \mapsto \ \xi - r \sum_{i|a_i<0} D_i$$
in the study of $\cL$ in Section \ref{grrrr}. 
The lowest $r$ terms of 
$$c_g(-R\pi_*\cL) \ \ \text{and} \ \ c_g(-R\pi_*\cL^{\text{orb}})$$
are therefore %DZ: removed ``exactly''.
equal, and we can apply Corollary
\ref{Cor:pushforwardck2} 
 to calculate
 $$
 \text{Coeff}_{r^0}
\left[
r \epsilon_* c_g(-R\pi_*\cL^{\text{orb}})\cdot \left[ \oM_{g,n,\beta}(X) \right]^\vir
\right]
\, \in\, A_{*}(\oM_{g,n,\beta}(X))\, .$$

Corollary \ref{Cor:pushforwardck2} gives the coefficient of $r^{-1}$ of $\epsilon_* c_g(-R\pi_*\cL)$ in the Artin stack $\ooMM_{g,n}^{Z}$.
The answer is obtained by the
 $r=0$ restriction of the degree~$g$ part of 
\begin{multline*}
\sum_{\Gamma\in \G_{g,n}^{Z}}\,
\sum_{r{\text{-twist}}\, \tw} 
%w\in \mathsf{W}_{\Gamma,r}
 \frac{r^{-h^1(\Gamma)}}{|\Aut(\Gamma)|}
\;
j_{\Gamma*}\Bigg[
%\prod_{i=1}^n \exp\left(\frac12 a_i^2 \psi_i + a_i \xi_i \right)
\prod_{v \in \V(\Gamma)} \exp\left(-\frac12 \eta(v) \right)
\\ \hspace{+10pt}
\prod_{e=(h,h')\in \E(\Gamma)}
\frac{1-\exp\left(-\frac{w(h)w(h')}2(\psi_h+\psi_{h'})\right)}{\psi_h + \psi_{h'}} \Bigg]\, .
\end{multline*} 
By applying Lemma~\ref{Lem:pullback1},
 we can calculate the pull-back of the class in  $\oM_{g,n,\beta}(X)$.
After interpreting the twists as weights,
 we obtain the $r=0$ restriction of the degree~$g$ part of
\begin{multline*}
\sum_{
\substack{\Gamma\in \G_{g,n,\beta}(X) \\
w\in \mathsf{W}_{\Gamma,r}}
}
 \frac{r^{-h^1(\Gamma)}}{|\Aut(\Gamma)|}
\;
j_{\Gamma*}\Bigg[
\prod_{i=1}^n \exp\left(\frac12 a_i^2 \psi_i + a_i \xi_i \right)
\prod_{v \in \V(\Gamma)} \exp\left(-\frac12 \eta(v) \right)
\\ \hspace{+10pt}
\prod_{e=(h,h')\in \E(\Gamma)}
\frac{1-\exp\left(-\frac{w(h)w(h')}2(\psi_h+\psi_{h'})\right)}{\psi_h + \psi_{h'}} \Bigg]\, .
\end{multline*} 
The result is exactly the formula for the $X$-valued ${\mathsf{DR}}$-cycle 
claimed in Theorem~\ref{Thm:main}. 
\qed

\section{Applications}
\label{appp}

\subsection{A topological view}
Let $X$ be a nonsingular projective variety with a line bundle
$S \rightarrow X$, and let 
$$\PP(\cO_X \oplus S) \rightarrow X$$
be the canonically associated $\CP^1$-bundle over $X$
with $0$-divisor $D_0$ and $\infty$-divisor
$D_\infty$.

Localization with respect to the fiberwise $\mathbb{C}^*$-action 
 immediately leads to a calculation of the
Gromov-Witten theory of  $\PP(\cO_X \oplus S)$ in terms of
the Gromov-Witten theory of
$X$ and the class $$c_1(S)\in H_2(X,\mathbb{Z})\, .$$
In \cite{MauPan}, an effective procedure was given to compute
the Gromov-Witten {invariants} of the associated
rubber geometry and the three relative geometries
$$\PP(\cO_X \oplus S)/D_0 \, \ , \ \ \
\PP(\cO_X \oplus S)/D_\infty\, \ , \ \ \
\PP(\cO_X \oplus S)/D_0 \cup D_\infty$$
in terms of the Gromov-Witten theory of $X$ and the class $c_1(S)$.
The results may be  viewed as analogues in Gromov-Witten theory
of the Leray-Hirsch
Theorem.

A basic consequence of
the $X$-valued $\mathsf{DR}$-cycle formula
of Theorem \ref{Thm:main} is a much stronger result on
the level of Gromov-Witten classes (not invariants).

\begin{proposition}\label{prfp}
The $X$-valued $\mathsf{DR}$-cycle formula
calculates the push-forward to the moduli space
of maps to $X$ of the virtual fundamental classes of
the moduli spaces of stable maps to
\begin{equation}\label{frr992}
\PP(\cO_X \oplus S)/D_0 \, \ , \ \ \
\PP(\cO_X \oplus S)/D_\infty\, \ , \ \ \
\PP(\cO_X \oplus S)/D_0 \cup D_\infty
\end{equation} 
in terms of tautological classes and $c_1(S)\in A^1(X)$.
\end{proposition}

\paragraph{Proof.}
Theorem \ref{Thm:main} provides a formula for the push-forward
to the moduli space of maps to $X$ of the virtual fundamental classes
of moduli space of
maps to rubber
in terms of tautological classes and $$c_1(S)\in A^1(X)\, .$$
To apply Theorem \ref{Thm:main}, we localize
the equivariant virtual fundamental classes of the
moduli spaces of stable maps to  the
three
relative geometries \eqref{frr992}
with respect to the
fiberwise $\mathbb{C}^*$-action. The $\mathbb{C}^*$-fixed contributions
are either absolute or relative. The virtual localization
formula \cite{GrP} on the absolute side is already of the desired form.
On the relative side, after removing the cotangent line via the
rubber calculus \cite[Section 1.5]{MauPan}, the desired form
is provided by Theorem \ref{Thm:main}. \qed
\vspace{8pt}

While \cite{MauPan} provides an algorithm for calculating 
the Gromov-Witten invariants of the relative
geometries \eqref{frr992}, the complexity of the method is not 
practical for calculations{\footnote{The results, however, have been used
for theoretical purposes, see \cite{PanDon,PixP}.}}. On the other hand,
Theorem \ref{Thm:main} may be used effectively for calculation.

In case $X$ is a point, exact calculations using Pixton's formula were
presented in \cite[Section 3]{JPPZ} for Hodge classes and
Hodge integrals. In Section \ref{kk8822} below,
an application of Theorem \ref{Thm:main} is presented where
 $X$ is the resolution of the
surface $A_\ell$-singularity.

\subsection{Resolution of surface singularities}
\label{kk8822}

%\subsubsection{Overview}
%Maulik \cite{Maulik} has computed the Gromov-Witten invariants of the toric surfaces $X_p$ obtained by 
%resolving an $A_p$ singularity. He also computed basic rubber invariants over $X_p$
% with target~$X$. Here we re-derive the second result from the first using the formula for the DR-cycle.

\subsubsection{Gromov-Witten invariants of $X_\ell$} 

Maulik \cite{Maulik} computed the Gromov-Witten invariants of the toric 
surface $X_\ell$ obtained by 
resolving the surface $A_\ell$-singularity.
We briefly review the geometry of the problem and refer the reader to \cite{Maulik} for a more
detailed treatment.
% We use the terminology for a general simple singularity, because some of the results, though not the methods that use equivariant cohomology and localization, generalize to simple singularities with simply laced Lie algebras (i.e., those whose all roots have the same length).

The resolution of the $A_\ell$-singularity is a nonsingular quasi-projective
 surface $X_\ell$ with $\ell$ exceptional divisors. 
The intersection pairing of the divisors is given by
 the Cartan matrix $C$ of the Lie algebra $A_\ell$. 
For a simply laced Lie algebra, the Cartan matrix is given by 
$$C_{ii} = -2\, , \ \ \ C_{ij} = 1$$
 if vertices $i$ and $j$ in the Dynkin diagram are connected by an edge.

There is a $(\C^*)^2$-action on the resolution $X_\ell$ of the $A_\ell$-singularity
which leaves every exceptional divisor invariant. 
We denote by $t_1$ and $t_2$ the corresponding equivariant weights. 
Because $X_\ell$ is a holomorphic symplectic variety, 
the ordinary Gromov-Witten invariants of~$X_\ell$ vanish except in degree~0. 
Maulik computed the  {\em reduced} Gromov-Witten invariants which
correspond to $(t_1+t_2)$-coefficient of the $(\C^*)^2$-equivariant Gromov-Witten invariants of $X_\ell$. The $(t_1+t_2)$-coefficient
is the lowest nonvanishing coefficient.

\begin{theorem}[Maulik \cite{Maulik}] \label{mmm123}
  Let $\alpha$ be a root, let
  $\beta = d \alpha\in H_2(X_\ell,\mathbb{Z})$
be a nonzero curve class, 
and let
$$\omega_1, \dots, \omega_p \in H^2(X_\ell,\mathbb{Z})$$
be divisor classes.  
 Let $b_1, \dots, b_p \geq 0$ and $c_1, \dots, c_q > 0$ be integers
subject to the dimensional constraint
$$
\sum_{i=1}^p b_i + \sum_{j=1}^{q} c_j = g+q\, .
$$
Then, we have 
\begin{eqnarray*}
\left< 
\prod_{i=1}^p \tau_{b_i}(\omega_i) \; \prod_{j=1}^q \tau_{c_j}(1)
\right>^{\!\! X_\ell,\,  \red}_{\!\!g,p+q,\beta} \; &=& \quad \frac{(2g + p + q - 3)!}{(2g + p - 3)!} \; d^{2g+p-3}  \\
& &\cdot 
\prod_{i=1}^p
\frac{b_i!}{(2b_i + 1)!} \left(- \frac12 \right)^{b_i}
(\alpha, \omega_i) \\
& & \cdot
\prod_{j=1}^q 
\frac{(c_j - 1)!}{(2c_j - 1)!} \left(- \frac12 \right)^{c_j-1} .
\end{eqnarray*}
If $\beta$ is not a multiple of a root or if the dimensional constraint is not satisfied, then the invariant vanishes.
\end{theorem}

\subsubsection{The $\mathsf{DR}$-cycle}
The rank of the Cartan matrix for $A_\ell$ is $\ell$.
The dimension of the equivariant cohomology of $X_\ell$ is $\ell+1$,
and the Poincar\'e intersection form is given by
an extension of the Cartan matrix $C$ :
\begin{equation} \label{Eq:eta}
\eta = \left(
\begin{array}{c|c}
\frac1{(p+1) t_1 t_2} & \displaystyle 0\\
\hline
\displaystyle 0 & \displaystyle C
\end{array}
\right).
\end{equation}

Let $A = (a_1, \dots, a_n)$ be a vector of integers
satisfying $\sum_{i=1}^n a_i = 0$. Let $\beta \neq 0$.
We can compute the reduced rubber Gromov-Witten invariant
\begin{equation} \label{Eq:DRInt}
\left< 
%\prod_{i=1}^\ell \tau_{b_i}(1) \; 
\prod_{i=1}^n \tau_0(\omega_i)
\cdot \DR_{g,A,\beta}(X_\ell, \mathcal{O}_{X_\ell}) \right>^{\!\!\red}
\end{equation}
using the formula of Theorem \ref{Thm:main}
for the $X_\ell$-valued $\mathsf{DR}$-cycle.
The intersection number \eqref{Eq:DRInt} is, by definition,
the coefficient of $t_1 + t_2$ in the corresponding equivariant intersection
(the lowest nonvanishing coefficient).

\begin{lemma} \label{Lem:onevertex}
  The only $X_\ell$-valued stable graphs which
  contribute to the coefficient of $t_1+t_2$ are graphs with one vertex.
\end{lemma}

\paragraph{Proof.}
The $(\mathbb{C}^*)^2$-equivariant Gromov-Witten invariant associated to a vertex $v$ of the graph $\Gamma$ with $\beta(v)\neq 0$ is a polynomial divisible by $t_1+t_2$ since $X_\ell$ is holomorphic symplectic.
By formula \eqref{Eq:eta} for $\eta$,
an edge of $\Gamma$ contributes either a constant, if the markings at the half-edges are divisors, or a factor of $$(p+1)t_1t_2\, ,$$
if the markings at the half-edges are equal to 1.

For a vertex $v$ of $\Gamma$ with $\beta(v)=0$, a more careful study using
the $(\mathbb{C}^*)^2$-localization formula for the degree 0 Gromov-Witten
invariants of $X_\ell$ is required.
The $(\C^*)^2$-invariant locus in $X_\ell$ consists of $\ell+1$ points.
The tangent weights at the $k$th point are
\begin{equation}\label{dd337}
\alpha_k(t_1,t_2) = (p+2-k)t_1 -(k-1)t_2\, , \ \ \ 
-\alpha_{k+1}(t_1, t_2) = (k-p-1)t_i + k t_2\,
\end{equation}
and satisfy the equation
$$\alpha_k-\alpha_{k+1}= t_1+t_2\, .$$
For a divisor class $\omega$ of $X_\ell$, we denote by $\omega^{(k)}(t_1,t_2)$
the restriction of $\omega$ to the $k$th invariant point in
$(\mathbb{C}^*)^2$-equivariant cohomology.
The restriction  $\omega^{(k)}(t_1,t_2)$
is a linear combination of $\alpha_k$ and $\alpha_{k+1}$.

The $(\C^*)^2$-fixed locus of $\oM_{\g(v),\n(v),\beta(v)=0}(X_\ell)$
is the union of $\ell+1$ copies of $\oM_{\g(v),\n(v)}$ corresponding to constant maps to the $\ell+1$ invariant points.
The contribution of the $k$th copy to the
$(\mathbb{C}^*)^2$-equivariant integral
$$
\left< \prod_{i=1}^p \tau_{b_i}(\omega_i) \; \prod_{j=1}^q \tau_{c_j}(1)
\right>^{X_\ell}_{\g(v),\, \n(v)=p+q,\, \beta(v)=0}
$$
equals
\begin{equation} \label{Eq:rationalf}
-\frac{\prod\limits_{i=1}^p \omega_i^{(k)}}{\alpha_k \, \alpha_{k+1}}
\int\limits_{\oM_{\g(v),\n(v)=p+q}}
\!\!\!\!
\prod_{i=1}^p \psi_i^{b_i}
\prod_{j=1}^q \psi_{p+j}^{c_j} \cdot 
\Lambda^\vee(\alpha_k) \Lambda^\vee(-\alpha_{k+1})\, .
\end{equation}
Here, we use the notation 
$$
\Lambda^\vee(\alpha) = \alpha^{\g(v)} - \alpha^{\g(v)-1} \lambda_1 + \dots 
+ (-1)^{\g(v)} \lambda_{\g(v)}\, .
$$
The weights $\alpha_k$, $\alpha_{k+1}$, and $\omega_i^{(k)}$ are linear forms in $t_1$ and $t_2$. 

By formula \eqref{dd337},  the weights  $\alpha_k$ and $\alpha_{k+1}$ 
are never proportional to $t_1+t_2$. Thus,
the $(t_1+t_2)$-valuation of the rational function~\eqref{Eq:rationalf} is nonnegative. In other words, every vertex $v$ of $\Gamma$ of nonzero degree has
$(t_1+t_2)$-valuation at least one, while every edge and degree zero vertex has $(t_1+t_2)$-valuation at least zero. Since we are interested in the coefficient of $t_1+t_2$ of the result,
there can be only one vertex of nonzero degree, and we can restrict ourselves to the $(t_1+t_2)$-valuation zero part of every degree zero vertex contribution. 

We can extract the $(t_1+t_2)$-valuation zero part of \eqref{Eq:rationalf} 
by substituting 
$$t_1 = t\, , \ \ t_2 = -t\, .$$
 In particular, then $\alpha_k = -\alpha_{k+1}$ and hence, by Mumford's 
identity \cite{Mum} for Hodge classes, 
$$
\Lambda^\vee(\alpha_k) \Lambda^\vee(-\alpha_{k+1}) = (-1)^{\g(v)} \alpha_k^{2\g(v)}\, .
$$
Thus, the contribution \eqref{Eq:rationalf} simplifies to
\begin{equation}\label{22ww2}
(-1)^{\g(v)-1} \alpha_k^{2\g(v)-2} \prod\limits_{i=1}^p \omega_i^{(k)}
\int\limits_{\oM_{\g(v),\n(v)=p+q}}
\!\!\!\!\prod_{i=1}^p \psi_i^{b_i}
\prod_{j=1}^q \psi_{p+j}^{c_j}\, ,
\end{equation}
where $\alpha_k$ and $\omega_i^{(k)}$ are now linear forms in~$t$. 

The contribution \eqref{22ww2}
 is a Laurent monomial in~$t$ of degree~$2g-2+p$. 
In the total contribution of the graph~$\Gamma$, we will have a product of these monomials over the genus~$0$ vertices and also a
product of monomials $$-(p+1)t^2$$ over the edges carrying the class $1$ on both half-edges. 
We distribute the edge factor $t^2$ to the two adjacent vertices, one
factor of $t$ for each vertex. For the invariant \eqref{Eq:DRInt},
 all legs of $\Gamma$ carry divisors. Hence, 
every half-edge carrying the class 1 contributes a factor of $t$ to the monomial. 
Every degree~$0$ vertex therefore 
contributes a factor of 
$$t^{2\g(v)-2+p+q} = t^{2\g(v)-2+\n(v)}\, .$$
By the stability condition,  
for a degree~0 vertex, the integer $2\g(v)-2+\n(v)$ is positive. 
Every degree~0 vertex thus contributes a monomial factor of positive degree. 
A product of such factors can never have a constant term.
We conclude that there are no degree~0 vertices and, 
as claimed, only one vertex of nonzero degree.
\qed 
\vspace{8pt}

We are ready now to compute the reduced rubber Gromov-Witten invariant 
$$
\left< 
%\prod_{i=1}^\ell \tau_{b_i}(1) \; 
\prod_{i=1}^n \tau_0(\omega_i)
\cdot \DR_{g,A,\beta}(X_\ell,\mathcal{O}_{X_\ell}) \right>^{\!\!\red}\, .
$$
For the computation, we will need two identities.

\begin{lemma} \label{Lem:summ}
We have
$$
\sum_{i+j = n} \frac1{(2i+1)! (2j+1)!} = \frac{2^{2n+1}}{(2n+2)!}\, .
$$
\end{lemma}

\paragraph{Proof.} After multiplying the left side by $(2n+2)!$, we obtain 
the well-known sum of odd binomial coefficients. \qed
\vspace{8pt}

The Bernoulli number are defined by the following generating series:
$$\sum_{m=0}^\infty B_m \frac{t^m}{m!} = \frac{t}{e^t-1}\, .$$
We define the functions $\cS(t)$ and $\cg(t)$ by
\begin{align*}
\cS(t) &= \frac{\sin(t/2)}{t/2} =\sum_{b \geq 0} \frac{(-1)^b}{(2b+1)! \, 4^b}\,  t^{2b}\, ,\\
\cg(t) &= - \frac12 \sum_{c \geq 1} \frac{B_{2c}}{2c} 
\sum_{c'+c'' = c-1} \frac{(-1)^{c'}}{(2c'+1)! \, 4^{c'}} \frac{(-1)^{c''}}{(2c''+1)! \, 4^{c''}} \,
t^{2c}\\
&= \sum_{c \geq 1} (-1)^c \frac{B_{2c}}{2c} 
\frac{t^{2c}}{(2c)!}\, .
\end{align*}
The last equality follows from Lemma~\ref{Lem:summ}.

\begin{lemma} \label{Lem:Sisexpg}
We have $\cS(t) = e^{\cg(t)}$.
\end{lemma}

\paragraph{Proof.}
The verification is straightforward:
\begin{eqnarray*}
t \left(\log(\cS(t))\right)' &=& t \left( \log \frac{\sin(t/2)}{t/2} \right)' \\
& = &
t \cdot \frac{t/2}{\sin(t/2)} \cdot \left(\frac{\frac12 \cos(t/2)}{t/2} 
- \frac{\frac12 \sin(t/2)}{(t/2)^2} \right) \\
&=& \frac{t}2 \cot(t/2) - 1\\
& =& \frac{it}{e^{it}-1}+\frac{it}2 -1 \\
& = &
\sum_{b \geq 1} (-1)^b B_{2b} \frac{t^{2b}}{(2b)!} \\
&=& t \cg'(t)\, .
\end{eqnarray*}
 \qed
\vspace{8pt}

We now substitute the formula of Theorem \ref{Thm:main} for
the $\mathsf{DR}$-cycle in \eqref{Eq:DRInt}.
Since only graphs with a single vertex contribute, we obtain
the following sum over the number $k$ of loops:
\begin{multline*}
%d^{2g+n-3} \;
\sum_{k \geq 0} \frac1{2^k \, k!} 
\;\;
\sum_{
\substack{
b_1, \dots, b_n\\
c_1, \dots, c_k}
}
\;\;
\sum_{
\substack{
c_1'+c_1''=c_1-1\\
\dots\\
c_k' + c_k'' = c_k-1}
}
\;\;
\sum_{
\substack{
\mu_1,\nu_1\\
\dots\\
\mu_k,\nu_k}
}
\;\;
\prod_{j=1}^k \left(-\eta^{\mu_j \nu_j} \frac{B_{2c_j}}{2c_j} \right)
\;\;
 \hspace{13em} \ 
\;\;
\\
\ \hspace{2em}
\times
\prod_{i=1}^n \frac{(a_i^2/2)^{b_i}}{b_i!} \prod_{j=1}^k \frac{(1/2)^{c_j'}}{c_j'!} \frac{(1/2)^{c_j''}}{c_j''!} \; 
\left< 
\prod_{i=1}^n \tau_{b_i}(\omega_i)
\prod_{j=1}^k \tau_{c_j'}(\omega_{\mu_j}) \tau_{c_j''}(\omega_{\nu_j})
\right>^{\text{red}}_{\!\!g-k,\beta} \!\!\!\!\!\!\!.
\end{multline*}

In the above sum, $c_j'$ and $c_j''$ are the powers of the $\psi$-classes at the branches of the $j$th node. The indices $\mu$ and $\nu$ run over a basis 
of $(\C^*)^2$-equivariant divisor classes of $X_\ell$
and $\eta^{\mu\nu}$ is the inverse of the Poincar\'e
intersection form $\eta_{\mu\nu}$.

In fact, we can  restrict the range of $\mu_j$ and $\nu_j$ in the
diagonal splitting by one.
Indeed, $$\eta^{1,1} = (p+1)t_1t_2\, ,$$
while we are interested in the coefficient of $t_1 + t_2$.
We can therefore replace the inverse of $\eta$ by the inverse of the Cartan matrix $C$ and write
$(C^{-1})^{\mu\nu}$ instead of $\eta^{\mu\nu}$.

Next, we apply Maulik's formula of Theorem \ref{mmm123} to
the reduced Gromov-Witten invariants which appear and use our
definitions of the generating series $\cS$ and $\cg$,
We conclude that \eqref{Eq:DRInt} equals 
\begin{multline*}
d^{2g+n-3} \prod_{i=1}^n (\alpha, \omega_i)
 \hspace{23em} \  \\
\ \hspace{3em}
\times
\sum_{k \geq 0} \frac1{2^k \, k!} \;
[t^{2g}] \left(\prod_{i=1}^n \cS(a_i t) \right)
\left( 2\cg(t)  \sum_{\mu,\nu} (\alpha, \omega_\mu) (C^{-1})^{\mu \nu} (\alpha, \omega_\nu)  \right)^k.
\end{multline*}
Since $\sum_{\mu,\nu} (\alpha, \omega_\mu) (C^{-1})^{\mu \nu} (\alpha, \omega_\nu) = (\alpha, \alpha) = -2$,
the formula simplifies to 
$$
d^{2g+n-3} \prod_{i=1}^n (\alpha, \omega_i)
\cdot 
\sum_{k \geq 0} \frac1{2^k \, k!} \;
[t^{2g}] \left(\prod_{i=1}^n \cS(a_i t) \right)
\left(-4 \cg(t) \right)^k 
$$
$$
= d^{2g+n-3} \prod_{i=1}^n (\alpha, \omega_i)
\cdot 
[t^{2g}] \left(\prod_{i=1}^n \cS(a_i t) \right) \exp(-2 \cg(t)) 
$$
$$
= d^{2g+n-3} \prod_{i=1}^n (\alpha, \omega_i) \cdot 
[t^{2g}]\frac{\prod_{i=1}^n \cS(a_i t) }{ \cS(t)^2}\, .
$$
The final equality coincides{\footnote{The result is also stated in
\cite[Proposition~3.1]{Maulik} where
    the factor $\prod_{i=1}^n (\alpha, \omega_i)$ is forgotten.}}
with  \cite[Proposition~3.6]{Maulik}
except for the automorphism factors of the partitions (omitted
here since we have numbered our marked points). 

The
same calculation as above using the formula of  Theorem \ref{Thm:main}
for the $\mathsf{DR}$-cycle yields the following more general evaluation
for the rubber theory over $X_\ell$.

\begin{theorem} \label{cc234}
Let $S \rightarrow X_\ell$ be a $(\mathbb{C}^*)^2$-equivariant
  line bundle on $X_\ell$. 
 Let $\alpha$ be a root,
let $\beta=d\alpha\in H_2(X_\ell,\mathbb{Z})$ be a nonzero curve
  class, and let
  $$A=(a_1,\ldots, a_n)$$ be a vector
  of integers satisfying
$$\sum_{i=1}^n a_i = \int_\beta c_1(S)\,. $$
Then, we have the evaluation
\begin{equation*}
\left< 
%\prod_{i=1}^\ell \tau_{b_i}(1) \; 
\prod_{i=1}^n \tau_0(\omega_i)
\cdot \DR_{g,A,\beta}(X_\ell, S) \right>^{\!\!\red} =
d^{2g+n-3} \prod_{i=1}^n (\alpha, \omega_i) \cdot 
[t^{2g}]\frac{\prod_{i=1}^n \cS(a_i t) }{ \cS(t)^2}\, .
\end{equation*}
\end{theorem}

\paragraph{Proof.}
Since $S\rightarrow X$ is now not necessarily trivial, Theorem \ref{Thm:main}
has additional $\xi_i$ terms at the markings and
$\pi_*(\xi^2)$ terms at the vertices. 
Lemma \ref{Lem:onevertex} still holds since the changes
in the $\mathsf{DR}$-cycle formula due to the line bundle
$S$ play no
role in the argument. The $k$-loop summation
for 
$$\left< 
%\prod_{i=1}^\ell \tau_{b_i}(1) \; 
\prod_{i=1}^n \tau_0(\omega_i)
\cdot \DR_{g,A,\beta}(X_\ell, S) \right>^{\!\!\red}$$
is again 
\begin{multline*}
%d^{2g+n-3} \;
\sum_{k \geq 0} \frac1{2^k \, k!} 
\;\;
\sum_{
\substack{
b_1, \dots, b_n\\
c_1, \dots, c_k}
}
\;\;
\sum_{
\substack{
c_1'+c_1''=c_1-1\\
\dots\\
c_k' + c_k'' = c_k-1}
}
\;\;
\sum_{
\substack{
\mu_1,\nu_1\\
\dots\\
\mu_k,\nu_k}
}
\;\;
\prod_{j=1}^k \left(-\eta^{\mu_j \nu_j} \frac{B_{2c_j}}{2c_j} \right)
\;\;
 \hspace{13em} \ 
\;\;
\\
\ \hspace{2em}
\times
\prod_{i=1}^n \frac{(a_i^2/2)^{b_i}}{b_i!} \prod_{j=1}^k \frac{(1/2)^{c_j'}}{c_j'!} \frac{(1/2)^{c_j''}}{c_j''!} \; 
\left< 
\prod_{i=1}^n \tau_{b_i}(\omega_i)
\prod_{j=1}^k \tau_{c_j'}(\omega_{\mu_j}) \tau_{c_j''}(\omega_{\nu_j})
\right>^{\text{red}}_{\!\!g-k,\beta} \!\!\!\!\!\!\! .
\end{multline*}
The extra $\xi_i$ and 
$\pi_*(\xi^2)$ terms produce additional factors of the 
equivariant parameters (and hence do not affect the reduced
invariants). 
The evaluation of the $k$-loop formula is then just
as before. \qed

\subsubsection{Remarks}
Maulik's evaluation of the reduced rubber invariants \eqref{Eq:DRInt}
played a crucial role in establishing the GW/DT/PT correspondences
for toric 3-folds, see \cite{MOOP,PPDesc}. His calculation of \eqref{Eq:DRInt}
in the case of $A_1$ relied upon the evaluation of the stationary theory of 
$\CP^1$ in \cite{OP1,OP2}. 
Using Maulik's $A_1$ argument in the reverse direction, 
Theorem \ref{Thm:main} via Theorem  \ref{cc234} 
provides a completely new $\mathsf{DR}$ derivation of the stationary Gromov-Witten theory of
$\CP^1$.

Theorem \ref{cc234} is also new. The $\mathsf{DR}$-cycle for
the rubber 
$$\mathbb{P}(\mathcal{O}_{X_\ell}\oplus S) \rightarrow X_\ell$$ 
constructed from the line bundle
$S \rightarrow X_\ell$
had not been considered before. The $\mathsf{DR}$
perspective puts all the rubber theories over $X_\ell$ on
the same footing.

\subsection{The tautological ring of the moduli of stable maps}
Pixton's formula
for the standard $\mathsf{DR}$-cycle leads to relations in the
tautological ring of $\oM_{g,n}$ first conjectured by Pixton \cite{PixRel}
and later
proven by Clader and Janda \cite{cj}.
After defining an appropriate strata algebra and tautological ring
for the moduli space of stable maps $\oM_{g,n,\beta}(X)$,  Bae \cite{Bae}
uses the formula of Theorem \ref{Thm:main} for the $X$-valued $\mathsf{DR}$-cycle
to construct tautological relations in the Chow theory of
$\oM_{g,n,\beta}(X)$, a rich new direction of study.

\subsection{Universal Abel-Jacobi theory on the Picard stack}
The classifying space of the group $\mathbb{C}^*$ is $\mathbb{CP}^\infty$.
If we take the  $\mathbb{CP}^N$-valued $\mathsf{DR}$-cycle in the limit
$N\rightarrow \infty$, we may hope that the result, suitably interpreted, is 
a universal $\mathsf{DR}$-cycle on the moduli space of line
bundles on curves. The required universal Abel-Jacobi theory on the Picard stack is developed in
\cite{BHPSS}. The calculation of the universal $\mathsf{DR}$-cycle there
is both motivated by and dependent upon our calculation of $\mathbb{CP}^N$-valued $\mathsf{DR}$-cycles.
The circle of ideas, also using \cite{HolSch},  leads to the proof of the formulas for the loci of
holomorphic and meromorphic differentials in $\oM_{g,n}$ conjectured in the Appendix of \cite{FarP}.

\newcommand{\arXiv}[1]{\texttt{arXiv:#1}}

\vspace{+12 pt}

\noindent Department of Mathematics, Univ. of Michigan\\
%\noindent Departement Mathematik \\
%\noindent ETH Z\"urich \\
\noindent janda@umich.edu 

\vspace{+6pt}
\noindent Departement Mathematik, ETH Z\"urich \\
\noindent rahul@math.ethz.ch

\vspace{+6 pt}
\noindent
Department of Mathematics, MIT\\
apixton@mit.edu

\vspace{+6 pt}
\noindent
CNRS, Versailles Mathematics Laboratory\\
dimitri.zvonkine@uvsq.fr

\end{document}